\documentclass[nonatbib, review]{elsarticle}
\journal{European Journal of Operational Research}
\usepackage{hyperref}
\usepackage{amsmath}
\usepackage{alphalph,etoolbox}
\usepackage{lipsum}
\usepackage{algorithm}
\usepackage{algorithmicx}
\usepackage{multicol}
\usepackage[labelformat=simple]{subcaption} 
\usepackage{tikz}
\usetikzlibrary{shapes.geometric, arrows}

\makeatletter
\renewcommand{\p@subfigure}{\thefigure} 
\makeatother

\usepackage{etoolbox} 
\AtBeginEnvironment{longtable}{\small}
\AtBeginEnvironment{table}{\small}

\patchcmd{\subequations}{\alph{equation}}{\alphalph{\value{equation}}}{}{}
\usepackage{blindtext}
\usepackage{anysize}
\usepackage{lmodern}
\usepackage{amsmath,amssymb,dsfont,amsthm,upgreek}

\usepackage{graphicx,color,tikz,float,rotating,enumitem}
\usepackage[draft]{changes}
\usepackage{comment}
\hypersetup{%
  colorlinks=true,
  linkcolor={blue!66!black},
  linkbordercolor=red,
  citecolor={blue!66!black}
  }
\definecolor{indigo}{RGB}{75,0,130}
\usepackage{longtable}
\usepackage{caption}
\usepackage{algpseudocode}
\usepackage{algorithm}
\usepackage{amsmath, amsthm, amssymb, bbm}

\usetikzlibrary{shapes.geometric, arrows}
\usepackage[title]{appendix}

\newcommand{\ignore}[1]{}

\usepackage{multirow}

\usepackage{apacite}
\usepackage[symbol]{footmisc}

\begin{document}

\begin{frontmatter}



\title{\textbf{Equilibrium models to analyse the impact of different coordination schemes between Transmission System Operator and Distribution System Operators on market power in sequentially-cleared energy and ancillary services markets under load and renewable generation uncertainty}}


\author[1]{Giovanni Micheli\footnote{Corresponding author.
\newline \hspace*{5.5mm} \textit{Email addresses}: giovanni.micheli@unibg.it (Giovanni Micheli),
maria-teresa.vespucci@unibg.it (Maria Teresa Vespucci),
gianluigi.migliavacca@rse-web.it (Gianluigi Migliavacca),
dario.siface@rse-web.it (Dario Siface)
}
}
\author[1]{Maria Teresa Vespucci}
\author[2]{Gianluigi Migliavacca}
\author[2]{Dario Siface}

\affiliation[1]{organization={Department of Management, Information and Production Engineering, University of Bergamo},
            city={Bergamo},
            country={Italy}}
\affiliation[2]{organization={RSE S.p.A},
            city={Milano},
            country={Italy}}

\begin{abstract}
The current massive installation of distributed resources in electricity distribution systems is transforming these systems into active dispatching subjects. At the same time, the need to compensate for the intermittent generation of an increasing amount of renewable sources creates the need to acquire more ancillary services. Flexible resources in the distribution system could provide these services not only within the perimeter of the distribution network to which they are connected but also for the benefit of the transmission system. However, this requires Transmission System Operators (TSOs) and Distribution System Operators (DSOs) to coordinate their dispatching actions effectively. One critical aspect of this coordination is establishing a market architecture that limits market power.

This paper presents an innovative game-theoretic approach to compare different TSO-DSO coordination models for acquiring ancillary services from distribution resources. Several schemes are considered: some with coordinated market management by TSOs and DSOs, others with sequential or independent local markets. For each scheme, the dispatching problem is formulated as a two-stage stochastic sequential game, where the first stage is the day-ahead market and the second stage is the balancing market. Nash equilibrium solutions are obtained by iteratively solving the profit maximization problem of each market player.
Numerical tests on a CIGRE benchmark network show that coordination schemes enabling distribution resources to provide ancillary services to the transmission system can significantly increase system costs when congestion occurs in the transmission network.

\end{abstract}


\end{frontmatter}



 \newpage
 
\section{Introduction}
\label{Intro}

The ambitious decarbonisation targets set by the European Commission and the need to make Europe as independent as possible from gas and fossil fuels are driving energy systems towards an ever increasing penetration of renewable energy sources, characterised by a typically intermittent generation pattern. The growth of Distributed Energy Resources (DERs), most frequently connected to the distribution networks, leads to an increased demand for flexibility services from system operators \cite{saha2023impact}.
However, it also creates an additional opportunity for DERs owners to provide flexibility services to the electricity system \cite{eid2016managing}, \cite{gulotta2023opening}. 

In some European countries, DERs are currently not allowed to participate in electricity markets \cite{arosio2023participation}. 
In Italy, for example, the dispatching rules state that only programmable generation units connected to the transmission grid can participate in the ancillary services market, while non-programmable renewable power plants are excluded from providing flexibility services \cite{terna2012grid}. 
However, the Italian regulator (ARERA) has initiated a process to open up flexibility markets to DERs \cite{ARERA2017}, and regulators in many other countries around the world, particularly in the EU, are moving in the same direction. Many EU Member States have already updated their national legislation to better integrate DERs into energy markets, and several countries are conducting experimental projects to assess the viability of new market structures \cite{rancilio2022ancillary}.
As stated in the European Directive 2019/944 \cite{EU}, the necessary regulatory framework has to be established by the Member States. As a result, many European Member States have launched a series of pilot projects to evaluate the possibility for DERs to provide balancing, congestion management and tertiary reserve services in the service market, see for example \cite{migliavacca2017smartnet,stanley2019platforms,utrilla2020large,losa2021platone,migliavacca2021innovative}.
A major limitation of existing studies is that they have largely neglected the potential exercise of market power by producers, i.e. the application of strategies to maximise the profit of the market participant, which on the other hand can lead to significant cost increases for the system.
The ability of dominant players to exercise market power must necessarily be considered when reforming the services market to allow for the participation of DERs.
In fact, energy markets are oligopolies with a limited number of participants: for example, in a hypothetical market for local distribution services, the radial topology of the distribution system greatly reduces the number of entities that can be active in resolving congestion in a given branch of the system.
Further opportunities for the exercise of market power arise from the existence of multiple markets in a cascade: it is well known that sequential trading can significantly affect competition and create arbitrage opportunities for market participants \cite{allaz1993cournot}.
Our research contributes to the ongoing study of the integration of DERs into energy markets by proposing a novel procedure that can assist transmission system operators (TSOs), distribution system operators (DSOs) and regulators in deciding on the design of the future electricity system and related energy markets. Specifically, we consider the possibility for DERs to provide congestion management and balancing services, and develop mathematical programming models to help regulators assessing how different modes of DER participation in the services market affect allocative efficiency and the potential for market power exploitation.
The specific feature of the proposed approach is the ability to capture the oligopolistic competition between market players participating in multi-level energy markets by formulating the problem as a "game" in which each subject brings its own optimisation problem and by determining a Nash equilibrium, i.e. a solution from which no agent is willing to deviate unilaterally.



To the best of our knowledge, there are no other works in the literature that use game theoretic models to investigate different coordination schemes between TSOs and DSOs in the provision of ancillary services. 
However, there are several examples of game theoretic models that have been developed to study the behaviour of electricity companies in deregulated electricity markets.
For example, in \cite{garci2006electricity} electricity markets with local marginal prices are investigated, in which three types of agents are involved, namely generators, consumer companies and an independent system operator.
The authors develop a procedure to determine a near-equilibrium solution, i.e. a solution where the Nash equilibrium condition is weakened to allow for the possibility that a player may have a small incentive to change his strategy.
Authors in \cite{hobbs2007nash} propose a Nash-Cournot approach to model competition among electricity generators. Both references \cite{garci2006electricity} and \cite{hobbs2007nash} focus on competition between market producers in single energy markets and thus do not consider sequentially cleared markets.
Other works deal with two-stage electricity markets. For instance, reference \cite{holmberg2015comparison} develops a game theoretic model to analyze the imperfections introduced by arbitrage opportunities. In the proposed modeling framework, generators are considered as price-takers. Reference \cite{dijk2011effect} extends the analysis presented in \cite{holmberg2015comparison} to include both imperfect competition and arbitrage opportunities. However, the performed analysis is limited to two-node networks.

To represent two-stage games, some papers in the literature develop two-stage stochastic problems with equilibrium constraints.
Two-stage stochastic problems are particularly suitable for modeling two-stage markets, as they can capture the interdependencies between the outcomes of the two stages. Some contributions have been proposed in the literature. For example, authors in \cite{zhang2013two} and \cite{zhang2010two} develop two-stage stochastic models with equilibrium constraints to analyze the strategic behaviour of producers in forward and spot markets. 
In \cite{hesamzadeh2018simulation} authors evaluate the design of a zonal electricity market with imperfect competition. The two-stage game is first formulated as a two-stage stochastic problem with equilibrium constraints, and then reformulated as a mixed-integer bilinear program, whose solution is computationally demanding.
In \cite{sarfati2018increase}, authors formulate a two-stage game to analyze imperfect competition of generators in zonal electricity markets with a day-ahead and a real-time market. 
The two-stage game is first formulated as a two-stage stochastic problem with equilibrium constraints, and then reformulated as a mixed-integer linear program (MILP) solved by using Benders Decomposition. 

Following an approach similar to that of \cite{sarfati2018increase}, in this work we formulate the competitive game between market players that participate to the energy markets as a multi-leader common follower game.
In the proposed approach, each market player is a leader, deciding about the bidding strategies of the resources under its control with the objective of maximizing its own profit.
On the other hand, the market operator is the follower, collecting bid prices submitted by market players and clearing the energy markets with the aim of minimizing total system costs.
Thus, the competitive game can be formulated as a set of interrelated, bilevel problems, one for each market player, where the upper-level represents the bidding problem faced by the market player and the lower-level is the clearing problem faced by the market operator.
By replacing the common lower-level with the respective \textit{Karush-Kuhn-Tucker} (KKT) optimality conditions, we obtain a set of Mathematical Programs with Equilibrium Constraints (MPECs) that we solve by adopting an iterative procedure.
This paper is organized as follows. Section \ref{sec:model} formulates the competitive game between market players under different configurations of the services market as MPECs. Section \ref{sec:tech} describes the modeling techniques adopted to linearize the bilinear terms and the complementary conditions so as to obtain equivalent MILP formulations.
The iterative procedure to determine a Nash equilibrium is presented in Section \ref{sec:algo}.
Numerical results for a case study are discussed in Section \ref{sec:numeric}. Finally, conclusions are drawn in Section 
\ref{sec:concl}.


\color{black}

\bigskip
\section{Problem description and formulation}
\label{sec:model}

In this section, we formulate the bidding problem of strategic players participating in the Day-Ahead Market (DAM) and in the Ancillary Services Markets (ASMs), which are settled in subsequent stages.
The power transmission network ${\cal T}$ is modeled as a set of lines $l \in {\cal L}^{\cal T}$ connecting nodes $n \in {\cal N}^{\cal T}$.  
Let ${\cal U}^{\cal T}$ and ${\cal R}^{\cal T}$ denote the sets of programmable and non-programmable generators, respectively, connected to network ${\cal T}$.  
Power distribution networks ${\cal D}_k$, $1 \leq k \leq K$, are connected to network ${\cal T}$. 
The distribution network ${\cal D}_k$ is modeled as a set of lines $l \in {\cal L}^{{\cal D}_k}$ connecting nodes $n \in {\cal N}^{{\cal D}_k}$.
Let ${\cal U}^{{\cal D}_k}$ and ${\cal R}^{{\cal D}_k}$ denote the sets of programmable and non-programmable generators, respectively, connected to network ${{\cal D}_k}$.
We also define the sets
${\cal N} = {\cal N}^{\cal T} \cup {\cal N}^{\cal D}$ with ${\cal N}^{\cal D} = \bigcup_{k=1}^K {\cal N}^{{\cal D}_k}$,
${\cal U} = {\cal U}^{\cal T} \cup {\cal U}^{\cal D}$ with ${\cal U}^{\cal D} = \bigcup_{k=1}^K {\cal U}^{{\cal D}_k}$, 
${\cal R} = {\cal R}^{\cal T} \cup {\cal R}^{\cal D}$ with ${\cal R}^{\cal D} = \bigcup_{k=1}^K {\cal R}^{{\cal D}_k}$,
and
${\cal L} = {\cal L}^{\cal T} \cup {\cal L}^{\cal D}$ with ${\cal L}^{\cal D} = \bigcup_{k=1}^K {\cal L}^{{\cal D}_k}$.
Let ${\cal U}_n$ and ${\cal R}_n$ denote the sets of programmable and non-programmable power plants, respectively, connected to node $n \in {\cal N}$.
The parameter $\delta_n$, $0 \leq \delta_n \leq 1$, expresses the maximum fraction of the load at node $n \in {\cal N}$ that can be curtailed in real time. If $\delta_n > 0$, the load at node $n$  is a \textit{flexible load}.

%

For the hour under consideration, given the load forecast $D_n$, $n \in {\cal N}$, the renewable production forecast $W_r$, $r \in {\cal R}$, and the sale bids  $(G_u,b_u^{\cal U})$, $u \in {\cal U}$, the DAM operator determines the accepted quantities $g_u$, $u \in {\cal U}$, of the sale bids so as to satisfy the net load $\sum_{n \in {\cal N}} D_n - \sum_{r \in {\cal R}}W_r$
while maximising the social benefit $SB=-\sum_{u \in {\cal U}} b_u^{\cal U} \ g_u$.
The DAM is a pay-as-clear market, the bid price $b_u^{\cal U}$ is the minimum price requested by power plant $u \in {\cal U}$ and the clearing price $\lambda$ is the optimal value of the dual variable associated with the constraint to satisfy the load. 
We assume that power plants compete on price, therefore, the bid quantity is equal to the plant capacity $G_u$.
Finally, note that we consider a bus-bar DAM (as it is, for instance, in France, Germany, or Spain).


However, the net load to be supplied in real time may differ from the day-ahead forecast, since 
the actual load $\tilde{D}_n$ may differ from the forecast $D_n$, and the actual renewable production $\tilde{W}_r$ may differ from the forecast $W_r$. 
Furthermore, the production plans defined by the DAM may not be feasible because network constraints have not been taken into account.
These issues are addressed by clearing the ASM in real time with the aim of recruiting resources to resolve the network congestion that may result from the provisional day-ahead dispatch and the imbalances, due to the differences between the realised net load and renewable production and the day-ahead forecast. 
A positive imbalance indicates that the real time net load exceeds the scheduled generation in the DAM, requiring the provision of additional generation resources on the service market: it can be managed either by increasing the production of programmable power plants or by curtailing the load in real time up to the maximum share $\delta_n$, $n \in \mathcal{N}$. 
Conversely, a negative imbalance indicates an excess of generation that must be reduced in real time and can be resolved either by reducing the production of the programmable power plants dispatched in the DAM or by curtailing renewable generation. 
Ancillary services markets are pay-as-bid markets. As with DAM, we assume that power plants compete on price and offer all available flexibility. 

Alternative coordination schemes between transmission and distribution service markets have been proposed.  
In this paper we analyse the following:
(A) a common ASM for transmission and distribution, in which all flexibility resources connected either to the transmission network or to a distribution network participate; 
(B) separate ASMs for transmission and distribution, with only resources connected to each network operating in each market; 
(C) separate ASMs for transmission and distribution, where distribution resources not used in the respective networks can provide services to the transmission network. 

%
%
In coordination scheme A, given  
the up-regulation bids $(G_u-g_u, \ b^{\cal U,\uparrow}_u)$, $u \in {\cal U}$, 
the down-regulation bids $(g_u, \ b^{\cal U,\downarrow}_u)$, $u \in {\cal U}$, and the load curtailment bids $(\delta_n \tilde{D}_n, \ b^{\cal N,\downarrow}_n)$, $n \in {\cal N}$, 
the ASM Operator determines the accepted quantities $g^{\uparrow}_u$, $g^{\downarrow}_u$, $u \in  {\cal U}$, and $d^{\downarrow}_n$, $n \in  {\cal N}$, as well as the renewable production curtailment $w^{\downarrow}_r$, $r \in  {\cal R}$, so as to solve at minimum cost the system imbalance 
\begin{equation}
\Delta = \sum_{n \in  {\cal N}} (\tilde{D}_n - D_n) + \sum_{r \in  {\cal R}} (W_r - \tilde{W}_r)    
\nonumber
\end{equation}
and ensure that in all lines $l \in {\cal L}$ the flows are within the limits. 
Note that $b^{\cal U,\uparrow}_u$, $b^{\cal U,\downarrow}_u$ and $b^{\cal N,\downarrow}_n$ denote the prices of the corresponding bids. 

%
%
In coordination scheme B, the ASM Operator for the transmission network $\cal T$ determines
$g^{\uparrow}_u$, $g^{\downarrow}_u$, $u \in  {\cal U^{\cal T}}$, $d^{\downarrow}_n$, $n \in  {\cal N^{\cal T}}$, and $w^{\downarrow}_r$, $r \in  {\cal R^{\cal T}}$,
given  
the up-regulation bids $(G_u-g_u, \ b^{\cal U,\uparrow}_u)$, $u \in {\cal U^{\cal T}}$, 
the down-regulation bids $(g_u, \ b^{\cal U,\downarrow}_u)$, $u \in {\cal U^{\cal T}}$, and the load curtailment bids $(\delta_n \tilde{D}_n, \ b^{\cal N,\downarrow}_n)$, $n \in {\cal N^{\cal T}}$, 
so as to solve at minimum cost the imbalance 
\begin{equation}
\Delta^{\cal T} = \sum_{n \in  {\cal N^{\cal T}}} (\tilde{D}_n - D_n) + \sum_{r \in  {\cal R^{\cal T}}} (W_r - \tilde{W}_r)    
\nonumber
\end{equation}
and ensure that in all lines $l \in {\cal L^{\cal T}}$ the flows are within the limits.
Parallely, for each distribution network ${\cal D}_k$, the ASM Operator  determines
$g^{\uparrow}_u$, $g^{\downarrow}_u$, $u \in {\cal U}^{{\cal D}_k}$, $d^{\downarrow}_n$, $n \in {\cal N}^{{\cal D}_k}$, and $w^{\downarrow}_r$, $r \in {\cal R}^{{\cal D}_k}$,
given  
the up-regulation bids $(G_u-g_u, \ b^{\cal U,\uparrow}_u)$, $u \in {\cal U}^{{\cal D}_k}$, 
the down-regulation bids $(g_u, \ b^{\cal U,\downarrow}_u)$, $u \in {\cal U}^{{\cal D}_k}$, and the load curtailment bids $(\delta_n \tilde{D}_n, \ b^{\cal N,\downarrow}_n)$, $n \in {\cal N}^{{\cal D}_k}$, 
so as to solve at minimum cost the network congestion, while ensuring that 
the exchange between ${\cal D}_k$ and  ${\cal T}$ after the ASM clearing equals the exchange between ${\cal D}_k$ and  ${\cal T}$ resulting from the DAM clearing.

%
%
In coordination scheme C, as in B, for each distribution network ${\cal D}_k$ the ASM Operator determines
$g^{\cal D,\uparrow}_u$, 
$g^{\cal D,\downarrow}_u$, $u \in {\cal U}^{{\cal D}_k}$, 
$d^{\cal D,\downarrow}_n$, $n \in {\cal N}^{{\cal D}_k}$, and $w^{\cal D,\downarrow}_r$, $r \in {\cal R}^{{\cal D}_k}$,
given
the up-regulation bids $(G_u-g_u, \ b^{\cal U, \cal D,\uparrow}_u)$, $u \in {\cal U}^{{\cal D}_k}$, 
the down-regulation bids $(g_u, \ b^{\cal U, \cal D,\downarrow}_u)$, $u \in {\cal U}^{{\cal D}_k}$, and 
the load curtailment bids $(\delta_n \tilde{D}_n, \ b^{\cal N, \cal D,\downarrow}_n)$, $n \in {\cal N}^{{\cal D}_k}$,
so as to solve at minimum cost the network congestion, while ensuring that 
the exchange between ${\cal D}_k$ and  ${\cal T}$ after the ASM clearing equals the exchange between ${\cal D}_k$ and ${\cal T}$ resulting from the DAM clearing.
Once the ASM of network ${\cal D}_k$ has been cleared, the flexible resources connected to ${\cal D}_k$ submit to the ASM of network $\cal T$ 
the up-regulation bids $(G_u-g_u-g^{\cal D,\uparrow}_u+g^{\cal D,\downarrow}_u, \ b^{\cal U, \cal T, \uparrow}_u)$, $u \in {\cal U}^{{\cal D}_k}$, 
the down-regulation bids $(g_u+g^{\cal D,\uparrow}_u-g^{\cal D,\downarrow}_u, \ b^{\cal U, \cal T, \downarrow}_u)$, $u \in {\cal U}^{{\cal D}_k}$, and 
the load curtailment bids $(\delta_n \tilde{D}_n-d^{\cal D,\downarrow}_n, \ b^{\cal N,\cal T,\downarrow}_n)$, $n \in {\cal N}^{{\cal D}_k}$.
In addition to the bids submitted by the flexible resources connected to networks ${\cal D}_k$, $1 \leq k \leq K$, the ASM Operator of network $\cal T$ receives the bids submitted by the flexible resources connected to network $\cal T$ and  
determines 
$g^{\cal T,\uparrow}_u$, 
$g^{\cal T,\downarrow}_u$, $u \in {\cal U}$, 
$d^{\cal T,\downarrow}_n$, $n \in {\cal N}$, and $w^{\cal T,\downarrow}_r$, $r \in {\cal R}$,
so as to solve the system imbalance $\Delta^{\cal T}$ and ensure that in all lines $l \in {\cal L^{\cal T}}$ the flows are within the limits.

%
%
The resources available either to supply load or to provide balancing and congestion management services are managed by a small number of strategic market players called \textit{Aggregators}. 
Let ${\cal I} = \{1,..., I\}$ denote the set of Aggregators,  
let ${\cal U}_i$ and ${\cal N}_i$ denote the sets of programmable generating units and flexible loads managed by Aggregator $i \in \cal I$
and let 
${\cal U}^{\cal D}_i \subseteq {\cal U}_i$ and ${\cal N}^{\cal D}_i \subseteq{\cal N}_i$ denote the sets of programmable generating units and flexible loads managed by Aggregator $i \in \cal I$ and connected to a distribution network.
Each Aggregator has to define joint bidding strategies on DAM and ASM, i.e., 
determine the bid prices 
$b^{\cal U}_u$, 
$b^{\cal U, \uparrow}_u$ and 
$b^{\cal U, \downarrow}_u$, for $u \in {\cal U}_i$, and  
$b^{\cal N, \downarrow}_n$, for $n \in {\cal N}_i$, in coordination schemes A and B, 
or 
$b^{\cal U}_u$, 
$b^{\cal U, \cal T, \uparrow}_u$ and 
$b^{\cal U, \cal T, \downarrow}_u$, for $u \in {\cal U}_i$,
$b^{\cal N, \cal T, \downarrow}_n$, for $n \in {\cal N}_i$,
$b^{\cal U, \cal D, \uparrow}_u$ and 
$b^{\cal U, \cal D, \downarrow}_u$, for $u \in {\cal U}^{\cal D}_i$,
and $b^{\cal N, \cal D, \downarrow}_n$, for $n \in {\cal N}^{\cal D}_i$,
in coordination scheme C, 
so as to maximize profit.
We assume that Aggregator $i \in \cal I$ selects each bid price from a set of discrete values: for instance, the bid price $b^{\cal U}_u$ of the sale bid on DAM related to the programmable generation unit $u \in {\cal U}_i$ is selected from the set $\{ B^{\cal U}_{u,a}$, $a \in {\cal A}^{\cal U}_u \}$.
We also assume that Aggregator $i$ determines the bidding strategy taking into account the uncertainty of the realised loads and renewable production, which is modelled by a set $\cal S$ of scenarios. Thus, we represent the bidding problem of Aggregator $i$ by a two-stage stochastic bilevel problem. The bid prices of Aggregator $i$, as well as the accepted quantities and the clearing price determined by the DAM operator, are first-stage variables; the accepted quantities of the up-regulation, down-regulation and load curtailment bids, as well as the renewable production curtailments, are second-stage variables.  
Finally, we assume that Aggregator $i$ knows the pricing strategies of the competitors, as observed in the market.

\subsection{Coordination scheme A: two-stage architecture}

In coordination scheme A,
Aggregator $i$ determines the bid prices 
$b^{\cal U}_u$, $b^{\cal U,\uparrow}_u$ and $b^{\cal U,\downarrow}_u$, for $u \in {\cal U}_i$, and $b^{\cal N,\downarrow}_n$, for $n \in {\cal N}_i$, by solving the following two-stage stochastic bilevel model \eqref{MODEL:ASM_2}:
\begin{subequations}
\label{MODEL:ASM_2}
{\allowdisplaybreaks
%
%
\begin{flalign}
\max
\sum_{u \in {\cal U}_i} (\lambda - C_u) \ g_u +
\sum_{s \in {\cal S}} \sigma_s \bigg\{ 
&
\sum_{u \in {\cal U}_i} 
\Big[  (b^{\cal U,\uparrow}_u - C^{\uparrow}_{u}) \ g^{\uparrow}_{u,s} +
(C^{\downarrow}_u - b^{\cal U,\downarrow}_u)  \ g^{\downarrow}_{u,s} 
\Big] +
\sum_{n \in {\cal N}_i} (b^{\cal N,\downarrow}_n - \lambda)
\ d^{\downarrow}_{n,s} \bigg\}
\label{eq:1a}
\end{flalign}
%
%
%
\begin{alignat}{2}
\text{s.t. } \quad
%
%
& b^{\cal U}_u = \sum_{a \in {\cal A}^{\cal U}_u} B^{\cal U}_{u,a} \ x^{\cal U}_{u,a}, && u \in {\cal U}_i,
\label{eq:1b} 
\\
%
%
& x^{\cal U}_{u,a} \in \{0, 1\}, \quad a \in {\cal A}^{\cal U}_u, && u \in {\cal U}_i, 
\label{eqn:1c} 
\\ 
%
%
& \sum_{a \in {\cal A}^{\cal U}_u} x^{\cal U}_{u,a} = 1, && u \in {\cal U}_i, 
\label{eq:1d}  
\\[10pt]
%
%
& g_u \in \arg \min \sum_{u \in {\cal U}} b^{\cal U}_u \ g_u
\label{eq:1e} \\
%
%
& \text{s.t. } 
  0 \leq g_u \leq G_u, && u \in {\cal U} \quad [\nu_u \geq 0], 
\label{eqn:1f} 
\\
%
%
& \qquad \sum_{u \in {\cal U}} g_u = \sum_{n \in {\cal N}} D_n -\sum_{r \in {\cal R}} W_r  && \qquad \quad \, \, \, [\lambda], 
\label{eqn:1g} 
\\[10 pt]
%
%
& b^{\cal U,\uparrow}_u = \sum_{a \in {\cal A}^{\cal U, \uparrow}_u} B^{\cal U, \uparrow}_{u,a} \ x^{\cal U, \uparrow}_{u,a}, && u \in {\cal U}_i,  
\label{eqn:1h} 
\\
%
%
& x^{\cal U, \uparrow}_{u,a} \in \{0, 1\}, && a \in {\cal A}^{\cal U, \uparrow}_u, u \in {\cal U}_i,  
\label{eqn:1i} 
\\
%
%
& \sum_{a \in {\cal A}^{\cal U, \uparrow}_u} x^{\cal U, \uparrow}_{u,a} = 1, && u \in {\cal U}_i,  
\label{eqn:1j} 
\\[10pt]
%
%
& b^{\cal U,\downarrow}_u = \sum_{a \in {\cal A}^{\cal U, \downarrow}_u} B^{\cal U, \downarrow}_{u,a} \ x^{\cal U, \downarrow}_{u,a}, && u \in {\cal U}_i,  
\label{eqn:1k} 
\\
%
%
& x^{\cal U, \downarrow}_{u,a} \in \{0, 1\}, && a \in {\cal A}^{\cal U, \downarrow}_u, u \in {\cal U}_i, 
\label{eqn:1l} 
\\
%
%
& \sum_{a \in {\cal A}^{\cal U, \downarrow}_u} x^{\cal U, \downarrow}_{u,a} = 1, && u \in {\cal U}_i,
\label{eqn:1m} 
\\[10pt]
%
%
& b^{\cal N,\downarrow}_n = \sum_{a \in {\cal A}^{\cal N, \downarrow}_n} B^{\cal N, \downarrow}_{n,a} \ x^{\cal N, \downarrow}_{n,a}, && n \in {\cal N}_i,
\label{eqn:1n} 
\\
%
%
& x^{\cal N, \downarrow}_{n,a} \in \{0, 1\}, && a \in {\cal A}^{\cal N, \downarrow}_n, n \in {\cal N}_i,
\label{eqn:1o} 
\\
%
%
& \sum_{a \in {\cal A}^{\cal N, \downarrow}_n} x^{\cal N, \downarrow}_{n,a} = 1, && n \in {\cal N}_i,
\label{eqn:1p} 
\\[10pt]
%
%
& \forall s \in {\cal S} \, (g^{\uparrow}_{u,s},d^{\downarrow}_{n,s},g^{\downarrow}_{u,s},w^{\downarrow}_{r,s}) \in \arg\min 
\sum_{u \in {\cal U}}  \Big( b^{\cal U,\uparrow}_u \ g^{\uparrow}_{u,s} - b^{\cal U,\downarrow}_u \ g^{\downarrow}_{u,s} \Big)
&& \nonumber \\
& \qquad \qquad \qquad \qquad \qquad \quad \qquad \qquad \,
+ \sum_{n \in {\cal N}} b^{\cal N,\downarrow}_n \ d^{\downarrow}_{n,s},
\label{eq:1q} 
\\
%
%
& \text{s.t. } \,
 0 \leq g^{\uparrow}_{u,s} \leq G_u - g_u, && u \in {\cal U} \quad [\beta_{u,s} \geq 0],
\label{eqn:1r} 
\\
%
%
& \qquad 0 \leq g^{\downarrow}_{u,s} \leq g_u, && u \in {\cal U} \quad [\phi_{u,s} \geq 0],
\label{eqn:1s} 
\\
%
%
& \qquad 0 \leq d^{\downarrow}_{n,s} \leq \delta_n \ \tilde{D}_{n,s}, && n \in {\cal N} \, \, \, \, \, [\gamma_{n,s} \geq 0],
\label{eqn:1t}
\\
%
%
& \qquad 0 \leq w^{\downarrow}_{r,s} \leq \tilde{W}_{r,s}, && r \in {\cal R} \quad [\chi_{r,s} \geq 0],
\label{eqn:1u}
\\
%
%
& \qquad
 \sum_{u \in {\cal U}} g^{\uparrow}_{u,s} + \sum_{n \in {\cal N}} d^{\downarrow}_{n,s} 
- \sum_{u \in {\cal U}} g^{\downarrow}_{u,s} - \sum_{r \in {\cal R}} w^{\downarrow}_{r,s} = \Delta_{s} && \qquad  \quad \, \, \, \, \, [\alpha_s],
\label{eqn:1v} 
\\
%
%
& \qquad \sum_{n \in {\cal N}} H_{l,n} \big[ \sum_{u \in {\cal U}_n} (g_u + g^{\uparrow}_{u,s} - g^{\downarrow}_{u,s}) && \nonumber \\
& \qquad 
+ \sum_{r \in {\cal R}_n}( \tilde{W}_{r,s} - w^{\downarrow}_{r,s}) - ( \tilde{D}_{n,s} - d^{\downarrow}_{n,s}) \big] \leq \overline{F}_l,  && l \in {\cal L} \quad [\mu_{l,s} \geq 0].
\label{eqn:1w} 
\end{alignat}
}
\end{subequations}


%
%
The upper level represents the bid selection problem of Aggregator $i$ and consists of the objective function \eqref{eq:1a} and the constraints 
\eqref{eq:1b}$-$\eqref{eq:1d} and
\eqref{eqn:1h}$-$\eqref{eqn:1p}.
%
%
Constraints \eqref{eq:1b} and \eqref{eqn:1c} express $b^{\cal U}_u$ as the linear combination of candidate bid prices $B^{\cal U}_{u,a}$, $a \in {\cal A}^{\cal U}_u$, where the binary variables $x^{\cal U}_{u,a}, a \in {\cal A}^{\cal U}_u$, are the combination coefficients. 
Constraint \eqref{eq:1d} imposes that exactly one candidate price is selected. 

%
%
Constraints \eqref{eq:1e}$-$\eqref{eqn:1g} represent the lower level problem solved by the DAM Operator,
given the bids $(G_u, b^{\cal U}_u), u \in \cal U$, of all Aggregators $i \in \cal I$ offering the capacity of their programmable power plants. We assume that Aggregator $i$ knows, or forecasts, the bid prices $b^{\cal U}_u, u \in \mathcal{U} \setminus {\mathcal U}_i$, of the competitors.
The accepted quantities $g_u, u \in \cal U$, and the market clearing price $\lambda$ are determined so as to maximize the social benefit \eqref{eq:1e}, subject to the constraints \eqref{eqn:1f}, imposing that $g_u$ is non negative and not greater than the offered quantity $G_u$, and constraint \eqref{eqn:1g}, imposing that the sum of the accepted quantities is equal to the net load, i.e. the difference between the total load $\sum_{n \in {\cal N}} D_n$ and the total non-programmable generation $\sum_{r \in {\cal R}} W_r$. 
The market clearing price $\lambda$ is the optimal value of dual variable associated to the balance constraint \eqref{eqn:1g}. 
 
%
%
Constraints 
\eqref{eqn:1h}$-$\eqref{eqn:1j}, 
\eqref{eqn:1k}$-$\eqref{eqn:1m} and 
\eqref{eqn:1n}$-$\eqref{eqn:1o}
impose the selection of bid prices 
$b^{\cal U,\uparrow}_u$, $b^{\cal U,\downarrow}_u$ and $b^{\cal N,\downarrow}_n$, respectively,
from the sets of candidates 
$\{B^{\cal U,\uparrow}_{u,a}, a \in {\cal A}^{\cal U,\uparrow}_u\}$, 
$\{B^{\cal U,\downarrow}_{u,a}, a \in {\cal A}^{\cal U,\downarrow}_u\}$ and
$\{B^{\cal N,\downarrow}_{n,a}, a \in {\cal A}^{\cal N,\downarrow}_n\}$, respectively.  

%
%
Constraints \eqref{eqn:1r}$-$\eqref{eqn:1w} represent the lower level problem solved in scenario $s \in \mathcal{S}$ by the Operator of the ASM common to networks $\cal T$ and ${\cal D}_k$, $1 \leq k \leq K$, given the up-regulation bids $(G_u-g_u, \ b^{\cal U,\uparrow}_u)$, $u \in {\cal U}$, the down-regulation bids $(g_u, \ b^{\cal U,\downarrow}_u)$, $u \in {\cal U}$, and the load curtailment bids $(\delta_n \tilde{D}_n, \ b^{\cal N,\downarrow}_n)$, $n \in {\cal N}$, submitted by Aggregators $i \in \cal I$ offering all their flexibility.
%
%
The accepted quantities $g^{\uparrow}_{u,s}$, $g^{\downarrow}_{u,s}$, $u \in  {\cal U}$, and $d^{\downarrow}_{n,s}$, $n \in  {\cal N}$, as well as the renewable production curtailment $w^{\downarrow}_{r,s}$, $r \in  {\cal R}$, are determined so as to minimize the total regulation cost \eqref{eq:1q}, i.e. the sum of the costs for up-regulations, minus the revenues of down-regulations, plus the costs of load curtailment. 
Constraints \eqref{eqn:1r}$-$\eqref{eqn:1u} impose that the accepted quantities and the renewable production curtailment do not exceed the offered quantities. 
Constraint \eqref{eqn:1v} imposes that the imbalance $\Delta_s$ generated in the real-time in scenario $s$
\begin{equation}
\Delta_s = \sum_{n \in  {\cal N}} (\tilde{D}_{n,s} - D_n) + \sum_{r \in  {\cal R}} (W_r - \tilde{W}_{r,s}) 
\nonumber
\end{equation}
is solved by modifying either the programmable generation, or the non-programmable generation, or the curtailable loads.
Constraint \eqref{eqn:1w} enforces the transit limits on power lines.
Variable $\alpha_s$ is the unrestricted dual variable of the balance constraint \eqref{eqn:1v}.
Variables $\beta_{u,s}$, $\phi_{u,s}$, $\gamma_{n,s}$, $\chi_{r,s}$ and $\mu_{l,s}$ are the dual variables associated with the upper bound constraints \eqref{eqn:1r}$-$\eqref{eqn:1u} and \eqref{eqn:1w}, respectively. 

%
%
The aggregator $i$ aims to find a profit maximising solution. 
In the objective function \eqref{eq:1a}, the first term represents the profit on the DAM, where $C_u$ is the operating cost of generation unit ${u \in {\cal U}_i}$, and the second term represents the expected profit on the ASM, where $C^\uparrow_u$ and $C^\downarrow_u$ are the up- and down-regulation costs of  generation unit ${u \in {\cal U}_i}$, respectively.  

%
%
Since the optimization problems \eqref{eq:1e}$-$\eqref{eqn:1g} 
and \eqref{eq:1q}$-$\eqref{eqn:1w} are linear, the associated \textit{Karush-Kuhn-Tucker} (KKT) conditions are both necessary and sufficient. 
Therefore, the stochastic bilevel program \eqref{MODEL:ASM_2} can be equivalently formulated as a single-level stochastic program with complementarity constraints. 
The equivalent formulation of problem \eqref{MODEL:ASM_2} as a single-level stochastic program with complementarity constraints is reported in \ref{KKT}.

%
%
\subsection{Coordination scheme B: three-stage architecture with distribution not supporting transmission}

As in coordination scheme A, 
the bid prices to be determined by Aggregator $i$ are
$b^{\cal U}_u$, $b^{\cal U,\uparrow}_u$ and $b^{\cal U,\downarrow}_u$, for $u \in {\cal U}_i$, and $b^{\cal N,\downarrow}_n$, for $n \in {\cal N}_i$, so as to maximise profit.
However, in coordination scheme B the ASM is divided into several markets for distribution (local) services and a market for transmission services. 
Only resources connected to the respective network operate in each market.
The following two-stage stochastic bilevel model \eqref{MODEL:ASM_3.1} represents the profit maximising bidding problem of Aggregator $i$:
\begin{subequations}
\label{MODEL:ASM_3.1}
{\allowdisplaybreaks
%
%
\begin{flalign}
\max
\sum_{u \in {\cal U}_i} (\lambda-C_u) \ g_u +
\sum_{s \in {\cal S}} \sigma_s \bigg\{ 
&
\sum_{u \in {\cal U}_i} 
\Big[  (b^{\cal U,\uparrow}_u - C^{\uparrow}_{u}) \ g^{\uparrow}_{u,s} +
(C^{\downarrow}_u - b^{\cal U,\downarrow}_u)  \ g^{\downarrow}_{u,s} \Big] +
\sum_{n \in {\cal N}_i} (b^{\cal N,\downarrow}_n - \lambda)
\ d^{\downarrow}_{n,s} \bigg\} 
\label{eq:2a}
\end{flalign}
\begin{alignat}{2}
%
%
\text{s.t. } \quad & \text{\eqref{eq:1b}$-$\eqref{eq:1d}}, 
\label{eq:2b}
\\
%
%
\qquad & \text{\eqref{eq:1e}$-$\eqref{eqn:1g}}, 
\label{eq:2c}
\\
%
%
\qquad & \text{\eqref{eqn:1h}$-$\eqref{eqn:1p}}, 
\label{eq:2d}
\\[10pt]
%
%
& \forall s \in {\cal S}, 1 \leq k \leq K \quad 
 (g^{\uparrow}_{u,s},g^{\downarrow}_{u,s},u \in {\cal U}^{{\cal D}_k},
\ d^{\downarrow}_{n,s}, n \in {\cal N}^{{\cal D}_k},
\ w^{\downarrow}_{r,s}, r \in {\cal R}^{{\cal D}_k}) \in \nonumber 
\\
& \arg\min 
  \sum_{u \in {\cal U}^{{\cal D}_k}}  \Big( b^{\cal U,\uparrow}_u \ g^{\uparrow}_{u,s} - b^{\cal U,\downarrow}_u \ g^{\downarrow}_{u,s} \Big) + 
  \sum_{n \in {\cal N}^{{\cal D}_k}} b^{\cal N,\downarrow}_n \ d^{\downarrow}_{n,s},
\label{eq:2e} \\
%
%
& \text{s.t. } \,
0 \leq g^{\uparrow}_{u,s} \leq G_u - g_u, && u \in {\cal U}^{{\cal D}_k},
\label{eqn:2f} \\
%
%
& \qquad 0 \leq g^{\downarrow}_{u,s} \leq g_u, && u \in {\cal U}^{{\cal D}_k},
\label{eqn:2g} 
\\
%
%
& \qquad 0 \leq d^{\downarrow}_{n,s} \leq \delta_n \ \tilde{D}_{n,s}, && n \in {\cal N}^{{\cal D}_k},
\label{eqn:2h}
\\
%
%
& \qquad 0 \leq w^{\downarrow}_{r,s} \leq \tilde{W}_{r,s}, && r \in {\cal R}^{{\cal D}_k},
\label{eqn:2i}
\\
%
%
& \qquad
 \sum_{u \in {\cal U}^{{\cal D}_k}} g^{\uparrow}_{u,s} + \sum_{n \in {\cal N}^{{\cal D}_k}} d^{\downarrow}_{n,s} 
- \sum_{u \in {\cal U}^{{\cal D}_k}} g^{\downarrow}_{u,s} - \sum_{r \in {\cal R}^{{\cal D}_k}} w^{\downarrow}_{r,s} = \Delta^{{\cal D}_k}_{s},
\label{eqn:2j} 
\\[3pt]
%
%
& \qquad \sum_{n \in {\cal N}^{{\cal D}_k}} H_{l,n} \big[ \sum_{u \in {\cal U}_n} (g_u + g^{\uparrow}_{u,s} - g^{\downarrow}_{u,s}) 
+ \sum_{r \in {\cal R}_n}( \tilde{W}_{r,s} - w^{\downarrow}_{r,s}) + && \nonumber \\
& \qquad \qquad - ( \tilde{D}_{n,s} - d^{\downarrow}_{n,s}) \big] \leq \overline{F}_l, && l \in {\cal L}^{{\cal D}_k},
\label{eqn:2k} 
\\[10pt]
%
%
& \forall s \in {\cal S} \quad (g^{\uparrow}_{u,s},g^{\downarrow}_{u,s},u \in {\cal U}^{\cal T},
\ d^{\downarrow}_{n,s}, n \in {\cal N}^{\cal T},
\ w^{\downarrow}_{r,s}, r \in {\cal R}^{\cal T}) \in \nonumber \\
& \arg\min 
  \sum_{u \in {\cal U}^{\cal T}}  \Big( b^{\cal U,\uparrow}_u \ g^{\uparrow}_{u,s} - b^{\cal U,\downarrow}_u \ g^{\downarrow}_{u,s} \Big) + 
  \sum_{n \in {\cal N}^{\cal T}} b^{\cal N,\downarrow}_n \ d^{\downarrow}_{n,s},
\label{eq:2l} 
\\
%
%
& \text{s.t. }
 \, 0 \leq g^{\uparrow}_{u,s} \leq G_u - g_u, && u \in {\cal U}^{\cal T},
\label{eqn:2m} 
\\
%
%
& \qquad 0 \leq g^{\downarrow}_{u,s} \leq g_u, && u \in {\cal U}^{\cal T},
\label{eqn:2n} 
\\
%
%
& \qquad 0 \leq d^{\downarrow}_{n,s} \leq \delta_n \ \tilde{D}_{n,s}, && n \in {\cal N}^{\cal T},
\label{eqn:2o}
\\
%
%
& \qquad 0 \leq w^{\downarrow}_{r,s} \leq \tilde{W}_{r,s}, && r \in {\cal R}^{\cal T},
\label{eqn:2p} 
\\
%
%
& \qquad
 \sum_{u \in {\cal U}^{\cal T}} g^{\uparrow}_{u,s} + \sum_{n \in {\cal N}^{\cal T}} d^{\downarrow}_{n,s} 
- \sum_{u \in {\cal U}^{\cal T}} g^{\downarrow}_{u,s} - \sum_{r \in {\cal R}^{\cal T}} w^{\downarrow}_{r,s} = \Delta^{\cal T}_{s},
\label{eqn:2q}
\\
%
%
& \qquad \sum_{n \in {\cal N}} H_{l,n} \big[ \sum_{u \in {\cal U}_n} (g_u + g^{\uparrow}_{u,s} - g^{\downarrow}_{u,s}) 
+ \sum_{r \in {\cal R}_n}( \tilde{W}_{r,s} - w^{\downarrow}_{r,s}) - ( \tilde{D}_{n,s} - d^{\downarrow}_{n,s}) \big] \leq \overline{F}_l, \,\, \,&&  l \in {\cal L}^{\cal T}.
\label{eqn:2r} 
\end{alignat}
}
\end{subequations}
Model \eqref{MODEL:ASM_3.1} differs from model \eqref{MODEL:ASM_2} as constraints \eqref{eq:1q}$-$\eqref{eqn:1w} are replaced by constraints \eqref{eq:2e}$-$\eqref{eqn:2r}.
In particular, 
constraints \eqref{eq:2e}$-$\eqref{eqn:2k} represent the lower level problem solved in scenario $s \in \mathcal{S}$ by the Operator of the ASM of network ${\cal D}_k$, $1 \leq k \leq K$, 
given the up-regulation bids $(G_u-g_u, \ b^{\cal U,\uparrow}_u)$, $u \in {\cal U}^{{\cal D}_k}$, the down-regulation bids $(g_u, \ b^{\cal U,\downarrow}_u)$, $u \in {\cal U}^{{\cal D}_k}$, and the load curtailment bids $(\delta_n \tilde{D}_n, \ b^{\cal N,\downarrow}_n)$, $n \in {\cal N}^{{\cal D}_k}$, submitted by Aggregators $i \in \cal I$.
The accepted quantities 
$g^{\uparrow}_{u,s}$, $g^{\downarrow}_{u,s}$, $u \in {\cal U}^{{\cal D}_k}$, and $d^{\downarrow}_{n,s}$, $n \in {\cal N}^{{\cal D}_k}$, as well as the renewable production curtailment $w^{\downarrow}_{r,s}$, $r \in {\cal R}^{{\cal D}_k}$, are determined so as to minimize the total regulation cost \eqref{eq:2e}.
Constraints \eqref{eqn:2f}$-$\eqref{eqn:2i} impose that the accepted quantities and the renewable production curtailment do not exceed the offered quantities.
Constraint \eqref{eqn:2j} ensures that after the ASM clearing the exchange between the distribution network ${\cal D}_k$ and the transmission network ${\cal T}$  equals the exchange resulting from the DAM clearing, since
\begin{equation*}
\Delta^{{\cal D}_k}_{s} = \sum_{n \in {\cal N}^{{\cal D}_k}} (\tilde{D}_{n,s} - D_n) + \sum_{r \in {\cal R}^{{\cal D}_k}} (W_r - \tilde{W}_{r,s}) 
\nonumber 
\end{equation*} 
and
\begin{equation*}
\sum_{u \in {\cal U}^{{\cal D}_k}} ( g_u + g^{\uparrow}_{u,s} - g^{\downarrow}_{u,s}) + \sum_{r \in {\cal R}^{{\cal D}_k}} ( \tilde{W}_{r,s} - w_{r,s})  - \sum_{n \in {\cal N}^{{\cal D}_k}} (\tilde{D}_{n,s} - d^{\downarrow}_{n,s} ) = \sum_{u \in {\cal U}^{{\cal D}_k}} g_u  + \sum_{r \in {\cal R}^{{\cal D}_k}}  W_{r,s} - \sum_{n \in {\cal N}^{{\cal D}_k}}{D}_{n,s}.     
\end{equation*} 
Constraint \eqref{eqn:2k} enforces the transit limits on power distribution lines.

Constraints \eqref{eq:2l}$-$\eqref{eqn:2r} represent the lower level problem solved in scenario $s \in \mathcal{S}$ by the Operator of the ASM of network $\cal T$, given the up-regulation bids $(G_u-g_u, \ b^{\cal U,\uparrow}_u)$, $u \in {\cal U}^{\cal T}$, the down-regulation bids $(g_u, \ b^{\cal U,\downarrow}_u)$, $u \in {\cal U}^{\cal T}$, and the load curtailment bids $(\delta_n \tilde{D}_n, \ b^{\cal N,\downarrow}_n)$, $n \in {\cal N}^{\cal T}$, submitted by Aggregators $i \in \cal I$.
The accepted quantities $g^{\uparrow}_{u,s}$, 
$g^{\downarrow}_{u,s}$, $u \in {\cal U}^{\cal T}$, 
and $d^{\downarrow}_{n,s}$, 
$n \in {\cal N}^{\cal T}$, 
as well as the renewable production curtailment $w^{\downarrow}_{r,s}$, $r \in {\cal R}^{\cal T}$, are determined so as to minimize the total regulation cost \eqref{eq:2l}.
Constraints \eqref{eqn:2m}$-$\eqref{eqn:2p} impose that the accepted quantities and the renewable production curtailment do not exceed the offered quantities.
Constraint \eqref{eqn:2q} states that the imbalance of the transmission network $\cal T$ must be solved by means of resources connected only to $\cal T$. 
Constraint \eqref{eqn:2r} enforces the transit limits on power distribution lines.

The equivalent formulation of problem \eqref{MODEL:ASM_3.1} as a single-level stochastic program with complementarity constraints is reported in \ref{KKT}.

%
%
\subsection{Coordination scheme C: three-stage architecture with distribution supporting transmission}

%
%
As in scheme B, also in coordination  scheme C, the ASM is split into several distribution (local) services markets and a transmission services market. 
Each DSO manages a local ASM to resolve congestion and imbalance in its distribution network at minimum cost, using only resources connected to its distribution network.
By contrast, all flexible resources (both transmission and distribution) can provide flexibility services to the TSO in order to resolve imbalances and congestion in the transmission system. 
Thus, distribution can support transmission by offering its residual resources for balancing and congestion management in transmission.

Each Aggregator $i$ determines the bid prices
$b^{\cal U}_u$, 
$b^{\cal U, \cal T, \uparrow}_u$ and 
$b^{\cal U, \cal T, \downarrow}_u$, for $u \in {\cal U}_i$,
$b^{\cal N, \cal T, \downarrow}_n$, for $n \in {\cal N}_i$,
$b^{\cal U, \cal D, \uparrow}_u$ and 
$b^{\cal U, \cal D, \downarrow}_u$, for $u \in {\cal U}^{\cal D}_i$,
and $b^{\cal N, \cal D, \downarrow}_n$, for $n \in {\cal N}^{\cal D}_i$,
by solving the following two-stage stochastic bilevel model \eqref{MODEL:ASM_3.2}:
\begin{subequations}
\label{MODEL:ASM_3.2}
{\allowdisplaybreaks
%
%
\begin{flalign}
\max \! 
\sum_{u \in {\cal U}_i}
(\lambda-C_u) g_u +
& \sum_{s \in {\cal S}} \sigma_s \bigg\{ 
\sum_{u \in {\cal U}^{\cal D}_i} 
\Big[ (b^{\cal U,\cal D,\uparrow}_u - C^{\uparrow}_{u}) g^{\cal D,\uparrow}_{u,s} + (C^{\downarrow}_u - b^{\cal U,\cal D,\downarrow}_u)  g^{\cal D,\downarrow}_{u,s} \Big] + \sum_{n \in {\cal N}^{\cal D}_i} (b^{\cal N,\cal D,\downarrow}_n - \lambda) d^{\cal D,\downarrow}_{n,s} 
\nonumber 
\\
& + \sum_{u \in {\cal U}_i}  \Big[ (b^{\cal U,\cal T,\uparrow}_u - C^{\uparrow}_{u}) g^{\cal T,\uparrow}_{u,s} +
(C^{\downarrow}_u - b^{\cal U,\cal T,\downarrow}_u)  g^{\cal T,\downarrow}_{u,s} \Big] + \sum_{n \in {\cal N}_i} (b^{\cal N,\cal T,\downarrow}_n - \lambda)
d^{\cal T,\downarrow}_{n,s} \bigg\} 
\label{eq:3a}
\end{flalign}
%
%
\begin{alignat}{2}
%
%
\text{s.t. } \quad & \text{\eqref{eq:1b}$-$\eqref{eq:1d}}, && 
\label{eq:3b}
\\
%
%
\quad & \text{\eqref{eq:1e}$-$\eqref{eqn:1g}}, && 
\label{eq:3c}
\\[10pt]
%
%
& {b}^{\cal U,\cal D,\uparrow}_u = \sum_{a \in {\cal A}^{\cal U, \uparrow}_u} B^{\cal U, \uparrow}_{u,a} \ {x}^{\cal U, \cal D,\uparrow}_{u,a}, && u \in {\cal U}^{\cal D}_i,  \label{eqn:3d} 
\\
%
%
& {x}^{\cal U, \cal D,\uparrow}_{u,a} \in \{0, 1\}, && a \in {\cal A}^{\cal U, \uparrow}_u, u \in {\cal U}^{\cal D}_i,  
\label{eqn:3e} 
\\
%
%
& \sum_{a \in {\cal A}^{\cal U, \uparrow}_u} {x}^{\cal U, \cal D,\uparrow}_{u,a} = 1, && u \in {\cal U}^{\cal D}_i, 
\label{eqn:3f} 
\\[10pt]
%
%
& {b}^{\cal U,\cal D,\downarrow}_u = \sum_{a \in {\cal A}^{\cal U, \downarrow}_u} B^{\cal U, \downarrow}_{u,a} \ {x}^{\cal U,\cal D, \downarrow}_{u,a}, &&  u \in {\cal U}^{\cal D}_i,  
\label{eqn:3g} 
\\
%
%
& {x}^{\cal U,\cal D,\downarrow}_{u,a} \in \{0, 1\}, && a \in {\cal A}^{\cal U, \downarrow}_u, u \in {\cal U}^{\cal D}_i, 
\label{eqn:3h} 
\\
%
%
& \sum_{a \in {\cal A}^{\cal U, \downarrow}_u} {x}^{\cal U,\cal D, \downarrow}_{u,a} = 1, && u \in {\cal U}^{\cal D}_i, 
\label{eqn:3i} 
\\[10pt]
%
%
& {b}^{\cal N,\cal D,\downarrow}_n = \sum_{a \in {\cal A}^{\cal N, \downarrow}_n} B^{\cal N, \downarrow}_{n,a} \ {x}^{\cal N,\cal D, \downarrow}_{n,a}, && n \in {\cal N}^{\cal D}_i,
\label{eqn:3j} 
\\
%
%
& {x}^{\cal N, \cal D,\downarrow}_{n,a} \in \{0, 1\}, &&  a \in {\cal A}^{\cal N, \downarrow}_n, n \in {\cal N}^{\cal D}_i,
\label{eqn:3k} 
\\
%
%
& \sum_{a \in {\cal A}^{\cal N, \downarrow}_n} {x}^{\cal N,\cal D, \downarrow}_{n,a} = 1, && n \in {\cal N}^{\cal D}_i, 
\label{eqn:3l}  
\\[10pt]
%
%
& \forall s \in {\cal S}, 1 \leq k \leq K \quad (g^{\cal D,\uparrow}_{u,s},g^{\cal D,\downarrow}_{u,s},u \in {\cal U}^{{\cal D}_k},
\ d^{\cal D,\downarrow}_{n,s}, n \in {\cal N}^{{\cal D}_k},
\ w^{\cal D,\downarrow}_{r,s}, r \in {\cal R}^{{\cal D}_k} &&) \in \nonumber \\
& \arg\min 
  \sum_{u \in {\cal U}^{{\cal D}_k}}  \Big( b^{\cal U,\cal D,\uparrow}_u \ g^{\cal D,\uparrow}_{u,s} - b^{\cal U,\cal D,\downarrow}_u \ g^{\cal D,\downarrow}_{u,s} \Big) + 
  \sum_{n \in {\cal N}^{{\cal D}_k}} b^{\cal N,\cal D,\downarrow}_n \ d^{\cal D,\downarrow}_{n,s},
\label{eq:3m} 
\\
%
%
& \text{s.t. } \,
 0 \leq g^{\cal D,\uparrow}_{u,s} \leq G_u - g_u, && u \in {\cal U}^{{\cal D}_k},
\label{eqn:3n} 
\\
%
%
& \qquad 0 \leq g^{\cal D,\downarrow}_{u,s} \leq g_u, && u \in {\cal U}^{{\cal D}_k},
\label{eqn:3o} \\
%
%
& \qquad 0 \leq d^{\cal D,\downarrow}_{n,s} \leq \delta_n \ \Tilde{D}_{n,s}, && n \in {\cal N}^{{\cal D}_k},
\label{eqn:3p}
\\
%
%
& \qquad 0 \leq w^{\cal D,\downarrow}_{r,s} \leq \Tilde{W}_{r,s}, && r \in {\cal R}^{{\cal D}_k},
\label{eqn:3q} 
\\
%
%
& \qquad \sum_{u \in {\cal U}^{{\cal D}_k}} g^{\cal D,\uparrow}_{u,s} + \sum_{n \in {\cal N}^{{\cal D}_k}} d^{\cal D,\downarrow}_{n,s} 
- \sum_{u \in {\cal U}^{{\cal D}_k}} g^{\cal D,\downarrow}_{u,s} - \sum_{r \in {\cal R}^{{\cal D}_k}} w^{\cal D,\downarrow}_{r,s} = \Delta^{{\cal D}_k}_{s},
\label{eqn:3r}
\\
%
%
& \qquad \! \sum_{n \in {\cal N}^{{\cal D}_k}} \! \! \! H_{l,n} \big[\! \! \sum_{u \in {\cal U}_n} (g_u + g^{\cal D,\uparrow}_{u,s} - g^{\cal D,\downarrow}_{u,s}) 
+ \! \! \sum_{r \in {\cal R}_n}( \Tilde{W}_{r,s} - w^{\cal D,\downarrow}_{r,s}) + && \nonumber \\
& \qquad - ( \Tilde{D}_{n,s} - d^{\cal D,\downarrow}_{n,s}) \big] \! \! \leq \overline{F}_l, && l \in {\cal L}^{{\cal D}_k},
\label{eqn:3s} 
\\[10pt]
%
%
& b^{\cal U,\cal T,\uparrow}_u = \sum_{a \in {\cal A}^{\cal U, \uparrow}_u} B^{\cal U, \uparrow}_{u,a} \ x^{\cal U,\cal T,\uparrow}_{u,a}, && u \in {\cal U}_i,  
\label{eqn:3t}
\\
%
%
& x^{\cal U, \cal T,\uparrow}_{u,a} \in \{0, 1\}, && a \in {\cal A}^{\cal U, \uparrow}_u, u \in {\cal U}_i,  
\label{eqn:3u} 
\\
%
%
& \sum_{a \in {\cal A}^{\cal U, \uparrow}_u} x^{\cal U,\cal T,\uparrow}_{u,a} = 1, && u \in {\cal U}_i,  
\label{eqn:3v} 
\\[10pt]
%
%
& b^{\cal U,\cal T,\downarrow}_u = \sum_{a \in {\cal A}^{\cal U, \downarrow}_u} B^{\cal U, \downarrow}_{u,a} \ x^{\cal U, \cal T,\downarrow}_{u,a}, && u \in {\cal U}_i,  
\label{eqn:3w} 
\\
%
%
& x^{\cal U, \cal T,\downarrow}_{u,a} \in \{0, 1\}, && a \in {\cal A}^{\cal U, \downarrow}_u, u \in {\cal U}_i, 
\label{eqn:3x} 
\\
%
%
& \sum_{a \in {\cal A}^{\cal U, \downarrow}_u} x^{\cal U, \cal T,\downarrow}_{u,a} = 1, && u \in {\cal U}_i, 
\label{eqn:3y} 
\\[10pt]
%
%
& b^{\cal N,\cal T,\downarrow}_n = \sum_{a \in {\cal A}^{\cal N, \downarrow}_n} B^{\cal N, \downarrow}_{n,a} \ x^{\cal N, \cal T,\downarrow}_{n,a}, && n \in {\cal N}_i,
\label{eqn:3z} 
\\
%
%
& x^{\cal N,\cal T,\downarrow}_{n,a} \in \{0, 1\}, && a \in {\cal A}^{\cal N, \downarrow}_n, n \in {\cal N}_i,
\label{eqn:3aa} 
\\
%
%
& \sum_{a \in {\cal A}^{\cal N, \downarrow}_n} x^{\cal N,\cal T, \downarrow}_{n,a} = 1, && n \in {\cal N}_i,
\label{eqn:3ab} 
\\[10pt]
%
%
& \forall s \in {\cal S} \quad (g^{\cal T,\uparrow}_{u,s},g^{\cal T,\downarrow}_{u,s},u \in {\cal U},
\ d^{\cal T,\downarrow}_{n,s}, n \in {\cal N},
\ w^{\cal T,\downarrow}_{r,s}, r \in {\cal R}) \in \nonumber \\
& \arg\min 
  \sum_{u \in {\cal U}}  \Big( b^{\cal U,\cal T,\uparrow}_u \ g^{\cal T,\uparrow}_{u,s} - b^{\cal U,\cal T,\downarrow}_u \ g^{\cal T,\downarrow}_{u,s} \Big) + 
  \sum_{n \in {\cal N}} b^{\cal N,\cal T,\downarrow}_n \ d^{\cal T,\downarrow}_{n,s},
\label{eq:3ac} 
\\
%
%
& \text{s.t. } \, 0 \leq g^{\cal T,\uparrow}_{u,s} \leq G_u - g_u - g^{\cal D,\uparrow}_{u,s} + g^{\cal D,\downarrow}_{u,s}, && u \in {\cal U}^{\cal D},
\label{eqn:3ad}
\\
%
%
& \qquad 0 \leq g^{\cal T,\uparrow}_{u,s} \leq G_u - g_u, && u \in {\cal U}^{\cal T},
\label{eqn:3ae}
\\
%
%
& \qquad 0 \leq g^{\cal T,\downarrow}_{u,s} \leq g_u + g^{\cal D,\uparrow}_{u,s} - g^{\cal D,\downarrow}_{u,s}, && u \in {\cal U}^{\cal D},
\label{eqn:3af} 
\\
%
%
& \qquad 0 \leq g^{\cal T,\downarrow}_{u,s} \leq g_u, && u \in {\cal U}^{\cal T},
\label{eqn:3ag} \\
%
%
& \qquad 0 \leq d^{\cal T,\downarrow}_{n,s} \leq \delta_n \ \Tilde{D}_{n,s} - d^{\cal D,\downarrow}_{n,s}, && n \in {\cal N}^{\cal D},
\label{eqn:3ah}
\\
%
%
& \qquad 0 \leq d^{\cal T,\downarrow}_{n,s} \leq \delta_n \ \Tilde{D}_{n,s}, && n \in {\cal N}^{\cal T},
\label{eqn:3ai}
\\
%
%
& \qquad 0 \leq w^{\cal T,\downarrow}_{r,s} \leq \Tilde{W}_{r,s}-w^{\cal D,\downarrow}_{r,s}, && r \in {\cal R}^{\cal D},
\label{eqn:3aj}
\\
%
%
& \qquad 0 \leq w^{\cal T,\downarrow}_{r,s} \leq \Tilde{W}_{r,s}, && r \in {\cal R}^{\cal T},
\label{eqn:3ak} 
\\[10pt]
%
%
& \qquad
 \sum_{u \in {\cal U}} g^{\cal T,\uparrow}_{u,s} + \sum_{n \in {\cal N}} d^{\cal T,\downarrow}_{n,s} 
- \sum_{u \in {\cal U}} g^{\cal T,\downarrow}_{u,s} - \sum_{r \in {\cal R}} w^{\cal T,\downarrow}_{r,s} = \Delta^{\cal T}_{s},
\label{eqn:3al} 
\\[10pt]
%
%
& \qquad 
\sum_{n \in {\cal N}^{\cal T}} H_{l,n} \big[ \sum_{u \in {\cal U}_n} (g_u + g^{\cal T,\uparrow}_{u,s} - g^{\cal T,\downarrow}_{u,s}) 
+ \sum_{r \in {\cal R}_n} (\Tilde{W}_{r,s} - w^{\cal T,\downarrow}_{r,s}) - ( \Tilde{D}_{n,s} + && \nonumber \\
& \qquad - d^{\cal T,\downarrow}_{n,s}) \big]  
+ \sum_{k=1}^K \sum_{n \in {\cal N}^{{\cal D}_k}} H_{l,n} \big[ \sum_{u \in {\cal U}_n} (g_u + g^{\cal T,\uparrow}_{u,s} - g^{\cal T,\downarrow}_{u,s} + {g}^{\cal D,\uparrow}_{u,s} -  {g}^{\cal D,\downarrow}_{u,s}) 
+ && \nonumber \\
& \qquad + \sum_{r \in {\cal R}_n}( \Tilde{W}_{r,s} 
- w^{\cal T,\downarrow}_{r,s} - {w}^{\cal D,\downarrow}_{r,s}) - ( \Tilde{D}_{n,s} - d^{\cal T,\downarrow}_{n,s} - {d}^{\cal D,\downarrow}_{n,s}) \big] \leq \overline{F}_l, && l \in {\cal L^T}.
\label{eqn:3am} 
\end{alignat}
}
\end{subequations}
Model \eqref{MODEL:ASM_3.2} differs from model \eqref{MODEL:ASM_2} as constraints \eqref{eqn:1h}$-$\eqref{eqn:1w} are replaced by constraints \eqref{eqn:3d}$-$\eqref{eqn:3am}.
%
%
Constraints 
\eqref{eqn:3d}$-$\eqref{eqn:3f}, 
\eqref{eqn:3g}$-$\eqref{eqn:3i} and 
\eqref{eqn:3j}$-$\eqref{eqn:3l}
impose the selection of bid prices 
$b^{\cal U, \cal D, \uparrow}_u$, 
$b^{\cal U, \cal D, \downarrow}_u$,
for $u \in {\cal U}^{\cal D}_i$,
and 
$b^{\cal N, \cal D, \downarrow}_n$, 
for $n \in {\cal N}^{\cal D}_i$,
respectively,
from the sets of candidates 
$\{B^{\cal U,\uparrow}_{u,a}, a \in {\cal A}^{\cal U,\uparrow}_u\}$, 
$\{B^{\cal U,\downarrow}_{u,a}, a \in {\cal A}^{\cal U,\downarrow}_u\}$ and
$\{B^{\cal N,\downarrow}_{n,a}, a \in {\cal A}^{\cal N,\downarrow}_n\}$, respectively. 

%
%
Constraints \eqref{eq:3m}$-$\eqref{eqn:3s} represent the lower level problem solved in scenario $s \in \mathcal{S}$ by DSO of network ${\cal D}_k$, $1 \leq k \leq K$, 
given 
the up-regulation bids $(G_u-g_u, \ b^{\cal U,\cal D,\uparrow}_u)$, $u \in {\cal U}^{{\cal D}_k}$, 
the down-regulation bids $(g_u, \ b^{\cal U,\cal D,\downarrow}_u)$, $u \in {\cal U}^{{\cal D}_k}$, and 
the load curtailment bids $(\delta_n \tilde{D}_n, \ b^{\cal N,\cal D,\downarrow}_n)$, $n \in {\cal N}^{{\cal D}_k}$, submitted by Aggregators $i \in \cal I$.
Network congestions and imbalances are resolved using only flexible resources connected to network ${\cal D}_k$.
The accepted quantities 
$g^{\cal D,\uparrow}_{u,s}$, 
$g^{\cal D,\downarrow}_{u,s}$, 
$u \in {\cal U}^{{\cal D}_k}$, and 
$d^{\cal D,\downarrow}_{n,s}$, $n \in {\cal N}^{{\cal D}_k}$, as well as the renewable production curtailment $w^{\cal D,\downarrow}_{r,s}$, $r \in {\cal R}^{{\cal D}_k}$, are determined so as to minimize the total regulation cost \eqref{eq:3m}.
Constraints \eqref{eqn:3n}$-$\eqref{eqn:3q} impose that the accepted quantities and the renewable production curtailment do not exceed the offered quantities.
Constraint \eqref{eqn:3r} ensures that after the ASM clearing the exchange between the distribution network ${\cal D}_k$ and the transmission network ${\cal T}$  equals the exchange resulting from the DAM clearing. 
Constraint \eqref{eqn:3s} enforces the transit limits on power distribution lines.

%
%
Constraints 
\eqref{eqn:3t}$-$\eqref{eqn:3v}, 
\eqref{eqn:3w}$-$\eqref{eqn:3y} and 
\eqref{eqn:3z}$-$\eqref{eqn:3ab}
impose the selection of bid prices 
$b^{\cal U, \cal T, \uparrow}_u$, 
$b^{\cal U, \cal T, \downarrow}_u$,
for $u \in {\cal U}_i$,
and 
$b^{\cal N, \cal D, \downarrow}_n$, 
for $n \in {\cal N}_i$,
respectively,
from the sets of candidates 
$\{B^{\cal U,\uparrow}_{u,a}, a \in {\cal A}^{\cal U,\uparrow}_u\}$, 
$\{B^{\cal U,\downarrow}_{u,a}, a \in {\cal A}^{\cal U,\downarrow}_u\}$ and
$\{B^{\cal N,\downarrow}_{n,a}, a \in {\cal A}^{\cal N,\downarrow}_n\}$, respectively. 

%
%
Constraints \eqref{eq:3ac}$-$\eqref{eqn:3am} represent the lower level problem solved in scenario $s \in \mathcal{S}$ by the TSO, given 
the up-regulation bids 
$(G_u-g_u-g^{\cal D,\uparrow}_{u,s}+g^{\cal D,\downarrow}_{u,s}, \ b^{\cal U, \cal T, \uparrow}_u)$, $u \in {\cal U}^{\cal D}$, 
and $(G_u-g_u, \ b^{\cal U, \cal T, \uparrow}_u)$, $u \in {\cal U}^{\cal T}$,
the down-regulation bids 
$(g_u+g^{\cal D,\uparrow}_{u,s}-g^{\cal D,\downarrow}_{u,s}, \ b^{\cal U, \cal T, \downarrow}_u)$, 
$u \in {\cal U}^{\cal D}$,
and $(g_u, \ b^{\cal U, \cal T, \downarrow}_u)$, $u \in {\cal U}^{\cal T}$, and
the load curtailment bids 
$(\delta_n \tilde{D}_{n,s}-d^{\cal D,\downarrow}_{n,s}, \ b^{\cal N,\cal T,\downarrow}_n)$, 
$n \in {\cal N}^{\cal D}$,
and
$(\delta_n \tilde{D}_{n,s}, \ b^{\cal N,\cal T,\downarrow}_n)$, 
$n \in {\cal N}^{\cal T}$.
Flexibility services to the transmission system can be provided either by transmission resources or distribution resources.
The accepted quantities 
$g^{\cal T,\uparrow}_{u,s}$, 
$g^{\cal T,\downarrow}_{u,s}$, $u \in {\cal U}$, 
$d^{\cal T,\downarrow}_{n,s}$, $n \in {\cal N}$, 
as well as the renewable production curtailment 
$w^{\cal T,\downarrow}_{r,s}$, $r \in {\cal R}$,
are determined so as to minimize the total regulation cost \eqref{eq:3ac}.
Constraints \eqref{eqn:3ad}$-$\eqref{eqn:3ak} impose that the accepted quantities and the renewable production curtailment do not exceed the offered quantities.
Constraint \eqref{eqn:3al} solves the imbalance of the transmission network $\cal T$ using flexibility resources connected both in transmission and in distribution. 
Constraint \eqref{eqn:3am} enforces the transit limits on power distribution lines.

The objective function \eqref{eq:3a} calculates the expected profit of Aggregator $i$, taking into account that flexible resources in the transmission network participate only in the market for transmission services, whereas flexible resources connected to the distribution network participate both in the market for distribution services and in the market for transmission services.

The equivalent formulation of problem \eqref{MODEL:ASM_3.2} as a single-level stochastic program with complementarity constraints is reported in \ref{KKT}.

\bigskip
\section{Modeling techniques to manage bilinear terms}
\label{sec:tech}

In Section \ref{sec:model}, the profit maximization problem of Aggregator $i$ in each coordination scheme is formulated as a mixed-integer stochastic bilevel problem. 
This formulation is then reformulated into a single-level optimization model by replacing the linear lower-level problems with their corresponding KKT conditions (see \ref{KKT}). 
However, solving the resulting optimization programs remains challenging due to the presence of bilinear terms in both the complementarity constraints and the objective function.
This section outlines the modeling techniques employed to linearize these bilinear terms, enabling the derivation of MILP formulations for the aggregator problems.

Regarding complementarity constraints, let us consider the primal constraint-dual variable pair:
$\sum_{j=1}^n a_{i,j}x_j \leq b_i \ : \ y_i \geq 0$. 
For the primal inequality constraint, we introduce the slack variable $s_i= b_i - \sum_{j=1}^n a_{i,j}x_j \geq 0$.
The complementarity condition, which requires the product of the slack variable and the dual variable to be zero (\(s_i y_i = 0\)), can be equivalently expressed by enforcing that \(s_i\) and \(y_i\) are Special Ordered Sets of type 1 (SOS1) variables.

With regard to the bilinear terms in the profit functions, it is possible to divide these terms into three categories: 
(i) products between ASM bid prices and accepted quantities on the ASM, 
(ii) products between DAM clearing price and accepted quantities on the DAM, and 
(iii) products between DAM clearing price and accepted quantities for curtailment bids on the ASM. 
The first category of bilinear terms includes products between quantities dispatched in the services market and prices of bids submitted in the same markets. 
Consider, for instance, coordination scheme A. 
The corresponding bilinear terms, see equation \eqref{eq:1a}, are 
$b^{\cal U,\uparrow}_u \ g^{\uparrow}_{u,s}$, 
$b^{\cal U,\downarrow}_u \ g^{\downarrow}_{u,s}$ and 
$b^{\cal N,\downarrow}_n \ d^{\downarrow}_{n,s}$. 
Applying the discretization constraints \eqref{eqn:1h}, \eqref{eqn:1k} and \eqref{eqn:1n}, 
these relations can be expressed as the product between continuous and binary variables. 
For instance, let us consider the upward regulation cost of unit $u \in {\cal U}$ in scenario $s \in {\cal S}$  
\begin{align}
b^{\cal U,\uparrow}_u \ g^{\uparrow}_{u,s} =  \sum_{a \in {\cal A}^{\cal U, \uparrow}_u} B^{\cal U, \uparrow}_{u,a} \ x^{\cal U, \uparrow}_{u,a} \ g^{\uparrow}_{u,s}, \quad u \in {\cal U}, s \in {\cal S}.
\end{align}
Applying McCormick's reformulation \cite{mccormick1976computability}, the bilinear terms $x^{\cal U, \uparrow}_{u,a} \ g^{\uparrow}_{u,s}$, which are in the form of the product between binary and continuous variables, can be linearized by introducing a new set of continuous variables $XG_{u,a,s}^{\uparrow}, u \in {\cal U}, a \in {\cal A}_u^{\cal U,\uparrow}, s \in {\cal S}$ and imposing the following additional constraints:
\begin{subequations}\label{eq:MC}
\begin{alignat}{2}
& g^{\uparrow}_{u,s} + G_u (x^{\cal U, \uparrow}_{u,a}-1) \leq XG_{u,a,s}^{\uparrow} \leq g^{\uparrow}_{u,s}, \qquad && u \in {\cal U}, a \in {\cal A}_u^{\cal U,\uparrow}, s \in {\cal S}, \label{eq:MC1} \\
& 0 \leq XG_{u,a,s}^{\uparrow} \leq G_u x^{\cal U, \uparrow}_{u,a}, \quad && u \in {\cal U}, a \in {\cal A}_u^{\cal U,\uparrow}, s \in {\cal S}. \label{eq:MC2}
\end{alignat}
\end{subequations}
Constraints \eqref{eq:MC} work as follows. When $x^{\cal U, \uparrow}_{u,a}=1$, constraint \eqref{eq:MC1} forces the equality between variables $XG_{u,a,s}^{\uparrow}$ and $g^{\uparrow}_{u,s}$, while constraint \eqref{eq:MC2} is redundant. On the other hand, for $x^{\cal U, \uparrow}_{u,a}=0$ constraint \eqref{eq:MC2} forces the variable $XG_{u,a,s}^{\uparrow}$ to take a null value, while constraint \eqref{eq:MC1} is redundant. 
The same considerations apply to the products  
$b^{\cal U,\downarrow}_u \ g^{\downarrow}_{u,s}$ and 
$b^{\cal N,\downarrow}_n \ d^{\downarrow}_{n,s}$. 

The second type of bilinear terms refers to the products between the quantities dispatched in the energy market $g_u$ and the market clearing price $\lambda$. 
Such bilinearities can be resolved using KKT conditions. 
For instance, let us consider coordination scheme A again: from complementarity conditions 
we have 
$g_u (b_u^{\cal U} + \nu_u - \lambda) = 0$ and thus $g_u \lambda = g_u b_u^{\cal U} + g_u \nu_u$. 
Since complementarity conditions also 
impose that $g_u \nu_u = G_u \nu_u$, the bilinear term in the objective function \eqref{eq:1a} can be formulated as
$g_u \lambda = g_u b_u^{\cal U} + G_u\nu_u, \ u \in {\cal U}_i$. 
The only bilinear term left (i.e., $g_u b_u^{\cal U}$) can be expressed as a product between continuous variables $g_u$ and binary variables $x_{u,a}^{\cal U}$ by applying discretisation constraint \eqref{eq:1b} and thus can be easily linearized by introducing a new set of auxiliary variables $XG_{u,a}$ and imposing McCormick constraints \eqref{eq:MC} as described above. 


The third category of bilinear terms includes products between the DAM clearing price $\lambda$ and the load curtailment on the ASM. 
In coordination scheme A, the corresponding bilinear term is $\lambda \ d^{\downarrow}_{n,s}$. 
Defining the clearing price as the price of the last accepted bid on the DAM, we can discretize the clearing price $\lambda$ by introducing the auxiliary binary variables $y_{u,a}, \ u \in {\cal U}, \ a \in {\cal A}_u^{\cal U}$ and imposing the following equations:
\begin{subequations}
\begin{align}
& \lambda = \sum_{u \in {\cal U}} 
\sum_{a \in {\cal A}_u^{\cal U}} B_{u,a}^{\cal U}  \ y_{u,a}, \label{eq:lambda_discr}\\
& \sum_{u \in {\cal U}} \sum_{a \in {\cal A}_u^{\cal U}} y_{u,a} = 1, \\
& y_{u,a} \in \{0,1\}.
\end{align}
\end{subequations}
Applying the discretization constraints \eqref{eq:lambda_discr}, 
the bilinear term $\lambda \ d^{\downarrow}_{n,s}$ can be expressed as the product between continuous and binary variables and linearized by imposing McCormick constraints \eqref{eq:MC} as described above.

\section{Procedure to determine a Nash equilibrium}
\label{sec:algo}

By applying the linearization techniques described in Section \ref{sec:tech}, we can formulate as a MILP model the profit maximization problem faced by aggregator $i$, hereafter denoted by ${\cal P}(i)$. 
Problem ${\cal P}(i)$ is non-convex and not differentiable due to the presence of binary variables, preventing the search for a Nash equilibrium through the formulation of the KKT conditions of all aggregators problems ${\cal P}(i), \ i \in \mathcal{I}$. 
In this section, we describe the iterative procedure we developed to determine a Nash equilibrium. 
In the proposed approach, we iteratively solve the profit maximization problem of each aggregator $i$, while having fixed the competitors bid prices, so as to determine the best response of aggregator $i$ to prices presented by the competing aggregators.
Once the decision problem ${\cal P}(i)$ of aggregator $i$ is solved, the prices for the bids submitted by resources $u \in {\cal U}_i, \ n \in {\cal N}_i$ are set to the optimal values determined by the solution of ${\cal P}(i)$ and the next aggregator $i+1$ is considered. 
The algorithm stops when the strategies of all aggregators do not change anymore: the solution found is by construction a Nash equilibrium as no aggregator is willing to deviate unilaterally. 
The proposed algorithm is described in detail by Algorithm \ref{algo:EQ}.
\begin{algorithm}
\caption{Equilibrium identification}
\label{algo:EQ}
\begin{multicols}{2}
\begin{algorithmic}[1]
\State Set the maximum number of iterations $K$
  \State Initialize bid prices of all flexible resources:
\State \quad \quad $b_u^{{\cal U},(0)} = 
\max_{a \in {\cal A}_u^{\cal U}} B_{u,a}, u \in {\cal U}_i$,
\State \quad \quad $b_u^{{\cal U},\uparrow,(0)} = 
\max_{a \in {\cal A}_u^{\cal U,\uparrow}} B_{u,a}^{\cal U,\uparrow}, u \in {\cal U}_i$,
\State \quad \quad $b_u^{{\cal U},\downarrow,(0)} = 
\min_{a \in {\cal A}_u^{\cal U,\downarrow}} B_{u,a}^{\cal U,\downarrow}, u \in {\cal U}_i$, 
\State \quad \quad $b_n^{{\cal N},\downarrow,(0)} = \max_{a \in {\cal A}_n^{\cal N,\downarrow}} B_{n,a}^{\cal N,\downarrow}, n \in {\cal N}_i$.
\State Initialize vector $\varepsilon_i=1, \ i \in {\cal I}$.
\For{$\kappa =1 : K$}
\State Retrieve bid prices from previous iteration:
\State \quad \quad $b_u^{{\cal U},(\kappa)} = b_u^{{\cal U},(\kappa-1)}, u \in {\cal U}$,
\State \quad \quad $b_u^{{\cal U},\uparrow,(\kappa)} = b_u^{{\cal U},\uparrow,(\kappa-1)}, u \in {\cal U}$,
\State \quad \quad $b_u^{{\cal U},\downarrow,(\kappa)} = b_u^{{\cal U},\downarrow,(\kappa-1)}, u \in {\cal U}$, 
\State \quad \quad $b_n^{{\cal N},\downarrow,(\kappa)} = b_n^{{\cal N},\downarrow,(\kappa-1)}, n \in {\cal N}$.
\For{$i = 1 : I$}
\State Fix competitors prices:
\State \quad \quad $b_u^{{\cal U}} = b_u^{{\cal U},(\kappa)}, u \in \mathcal{U} \setminus {\cal U}_i$,
\State \quad \quad $b_u^{{\cal U},\uparrow} = b_u^{{\cal U},\uparrow,(\kappa)}, u \in \mathcal{U} \setminus {\cal U}_i$,
\State \quad \quad$b_u^{{\cal U},\downarrow} = b_u^{{\cal U},\downarrow,(\kappa)}, u \in \mathcal{U} \setminus {\cal U}_i$,
\State \quad \quad $b_n^{{\cal N},\downarrow} = b_n^{{\cal N},\downarrow,(\kappa)}, n \in \mathcal{N} \setminus {\cal N}_i$.
\State Solve model ${\cal P}(i)$
\If{$b_u^{{\cal U},*} = b_u^{{\cal U},(\kappa)}, b_u^{{\cal U},\uparrow,*} = b_u^{{\cal U},\uparrow,(\kappa)}, b_u^{{\cal U},\downarrow,*} = b_u^{{\cal U},\downarrow,(\kappa)}, \ u \in {\cal U}_i$, and
$b_n^{{\cal N},\downarrow,*} = b_n^{{\cal N},\downarrow,(\kappa)}, \ n \in {\cal N}_i$}
  \State Set $\varepsilon_i = 0$
  \Else 
  \State Update bid prices to the optimal values determined by the model:
\State \quad \quad $b_u^{{\cal U},(\kappa)} = b_u^{{\cal U},*}, u \in {\cal U}_i$,
\State \quad \quad $b_u^{{\cal U},\uparrow,(\kappa)} = b_u^{{\cal U},\uparrow,*}, u \in {\cal U}_i$,
\State \quad \quad $b_u^{{\cal U},\downarrow,(\kappa)} = b_u^{{\cal U},\downarrow,*}, u \in {\cal U}_i$, 
\State \quad \quad $b_n^{{\cal N},\downarrow,(\kappa)} = b_n^{{\cal N},\downarrow,*}, n \in {\cal N}_i$.
\State Set $\varepsilon_i = 1$
  \EndIf
  \EndFor
  \If{$\sum_{i=1}^I \varepsilon_i = 0$}
  \State break
  \EndIf
  \EndFor
\end{algorithmic}
\end{multicols}
\end{algorithm}

The algorithm first requires the specification of a maximum number of iterations, denoted by $K$. In the initialization phase, maximum profit values are assigned to parameters indicating the prices in the current solution of the bids submitted by all flexible resources (i.e., maximum values for day-ahead, up-regulation and load curtailment bids, and minimum values for down-regulation and generation curtailment bids). The aim of this initialization is to increase the probability of finding the maximum profits equilibrium, in case of multiple equilibria. 
Each component of vector $\varepsilon_i$ controlling the algorithm convergence is set to 1. Finally, the initialization phase ends by setting the iteration counter $\kappa=1$.
The iterative step of the algorithm consists of an outer loop on the set of iterations and an inner loop on the set of aggregators. Specifically, at a given iteration $\kappa$ the first aggregator is considered and the respective decision problem is solved, while fixing decision variables associated with competitors' bid prices. If the aggregator's optimal strategy is unchanged from the previously determined strategy, the correspondent component of parameter $\varepsilon_i$ is assigned null value. On the other hand, if the aggregator's optimal response changes, parameters indicating the optimal bid prices are updated and the correspondent component of parameter $\varepsilon_i$ is assigned unit value, indicating the change in the aggregator's strategy.
When all aggregators are considered, we perform a convergence test by evaluating the vector $\varepsilon_i$. If all the components are zero, then the optimal response of each aggregator has not changed with respect to the previous iteration, indicating that an equilibrium solution has been reached, since no aggregator is willing to unilaterally deviate from this solution. Otherwise, it means that at least one aggregator has deviated in its optimal strategy from the solution at the previous iteration, indicating that an equilibrium solution has not been reached yet. In such a case, the algorithm continues by performing a new iteration.

\section{Numerical experiments}
\label{sec:numeric}

In this section, we present an extensive set of computational results to validate the proposed methodology.
The section is organized as follows.
Section \ref{sec:data} outlines the main features of the case study, while Section \ref{sec:results} presents the outputs of the analysis.

\subsection{The data}
\label{sec:data}

To provide the reader with a clear understanding of how different ASM configurations impact dispatch, this section presents results from a simplified system based on CIGRE benchmark networks \cite{strunztf}. 
Specifically, the transmission system comprises 12 nodes, with nodes 4, 5, and 6 each connected to distinct distribution networks. 
All distribution networks share an identical topology, comprising 13 lines and 14 nodes each.
The transmission network topology is depicted in Fig. \ref{fig:Tsystem}, while the common topology of the distribution networks is illustrated in Fig. \ref{fig:Dsystem}.
Despite its reduced complexity, this test case effectively highlights notable differences in the potential for market power exercise across various market configurations, emphasizing the critical role of the proposed models in supporting market design activities.

\begin{figure}
    \centering
    \resizebox{0.65\textwidth}{0.4\textwidth}{%
    \begin{tikzpicture}[thick]
    \draw[thick] (2, 0) -- (3.5, 0); 
    \node at (1.8,0) {9} ;
    \draw[thick] (4.5, 0) -- (6, 0); 
    \node at (4.3,0) {10} ;
    \draw[thick] (12, 0) -- (13.5, 0); 
    \node at (11.8,0) {11} ;

    \draw[thick] (2, -2) -- (3.5, -2); 
    \node at (1.8,-2) {1} ;
    \draw[thick] (4.5, -2) -- (6, -2); 
    \node at (4.3,-2) {2} ;
    \draw[thick] (7, -2) -- (8.5, -2); 
    \node at (6.8,-2) {5} ;
    \draw[thick] (9.5, -2) -- (11, -2); 
    \node at (9.3,-2) {4} ;
    \draw[thick] (12, -2) -- (13.5, -2); 
    \node at (11.8,-2) {3} ;

    \draw[thick] (2, -6) -- (3.5, -6); 
    \node at (1.8,-6) {7} ;
    \draw[thick] (12, -6) -- (13.5, -6); 
    \node at (11.8,-6) {8} ;

    \draw[thick] (5.75, -4) -- (7.25, -4); 
    \node at (5.55,-4) {12} ;
    \draw[thick] (5.75, -5) -- (7.25, -5); 
    \node at (5.55,-5) {6} ;

    \node[draw, circle,scale=0.8] (u1) at (2.5, 0.5) {U4}; %
    \node[draw, circle,scale=0.8] (u2) at (5, 0.5) {U1}; %
    \node[draw, circle,scale=0.8] (u3) at (12.5, 0.5) {U2}; %
    \node[draw, circle,scale=0.8] (u4) at (6.5, -3.5) {U3}; %

    \node[draw, circle,scale=0.8] (r1) at (2.75, -1.5) {R1}; %
    \node[draw, circle,scale=0.8] (r2) at (2.75, -5.5) {R2}; %
    \node[draw, circle,scale=0.8] (r3) at (12.5, -5.5) {R3}; %

    \draw (u1) -- (2.5, 0); 
    \draw (u2) -- (5, 0); 
    \draw (u3) -- (12.5, 0); 
    \draw (u4) -- (6.5, -4);

    \draw (r1) -- (2.75, -2);
    \draw (r2) -- (2.75, -6);
    \draw (r3) -- (12.5, -6);
    
    \draw[dashed, thick] (2.30, 0) -- (2.30, -2); 
    \draw[dashed, thick] (4.8, 0) -- (4.8, -2);
    \draw[dashed, thick] (13.2, 0) -- (13.2, -2);

    \draw[dashed, thick] (3.2, -2) -- (3.2, -2.3); 
    \draw[dashed, thick] (3.2, -2.3) -- (4.8, -2.3); 
    \draw[dashed, thick] (4.8, -2) -- (4.8, -2.3); 

    \draw[dashed, thick] (5.7, -2) -- (5.7, -2.3); 
    \draw[dashed, thick] (5.7, -2.3) -- (7.3, -2.3); 
    \draw[dashed, thick] (7.3, -2) -- (7.3, -2.3); 

    \draw[dashed, thick] (8.2, -2) -- (8.2, -2.3); 
    \draw[dashed, thick] (8.2, -2.3) -- (9.8, -2.3); 
    \draw[dashed, thick] (9.8, -2) -- (9.8, -2.3); 
    
    \draw[dashed, thick] (10.7, -2) -- (10.7, -2.3); 
    \draw[dashed, thick] (10.7, -2.3) -- (12.3, -2.3); 
    \draw[dashed, thick] (12.3, -2) -- (12.3, -2.3); 

    \draw[dashed, thick] (2.30, -2) -- (2.30, -6); 
    \draw[dashed, thick] (13.2, -2) -- (13.2, -6);

    \draw[dashed, thick] (2.75, -6) -- (2.75, -6.3); 
    \draw[dashed, thick] (2.75, -6.3) -- (12.75, -6.3); 
    \draw[dashed, thick] (12.75, -6) -- (12.75, -6.3); 

    \draw[dashed, thick] (2.75, -2) -- (2.75, -2.3);
    \draw[dashed, thick] (2.75, -2.3) -- (6, -5); 

    \draw[dashed, thick] (10.25, -2) -- (10.25, -2.3);
    \draw[dashed, thick] (10.25, -2.3) -- (7.0, -5); 

    \draw[dashed, thick] (6.5, -4) -- (6.5, -5);

    \draw[->,very thick] (3, 0) -- (3, 0.8); 
    \draw[->,very thick] (5.25, -2) -- (5.25, -1.2);
    \draw[->,very thick] (7.75, -2) -- (7.75, -1.2);
    \draw[->,very thick] (10.25, -2) -- (10.25, -1.2);
    \draw[->,very thick] (12.75, -2) -- (12.75, -1.2);
    \draw[->,very thick] (6.5, -5) -- (6.5, -5.8);

    \draw[->,very thick] (2.3, -6) -- (2.3, -6.8);
    \draw[->,very thick] (13.2, -6) -- (13.2, -6.8);
\end{tikzpicture}
}
    \caption{Topology of the transmission network. The transmission system consists of four programmable generators (U1, U2, U3, U4), three non-programmable generators (R1, R2, R3), and eight loads. Of these, five are flexible loads, located at nodes 2, 3, 4, 5, and 6.}
    \label{fig:Tsystem}
\end{figure}


\begin{figure}
    \centering
    \resizebox{0.5\textwidth}{0.6\textwidth}{%
    
\begin{tikzpicture}[thick]
    \draw[thick] (-0.75, 0) -- (0.75, 0); 
    \node at (-1,0) {13} ;
    \draw[thick] (-0.75, -1.5) -- (0.75, -1.5); 
    \node at (-1,-1.5) {14} ;
    \draw[thick] (-0.75, -3) -- (0.75, -3); 
    \node at (-1,-3) {15} ;
    \draw[thick] (-3, -4.5) -- (-1.5, -4.5); 
    \node at (-3.2,-4.5) {16} ;
    \draw[thick] (-3, -6) -- (-1.5, -6); 
    \node at (-3.2,-6) {17} ;
    \draw[thick] (-3, -7.5) -- (-1.5, -7.5); 
    \node at (-3.2,-7.5) {18} ;
    \draw[thick] (-3, -9) -- (-1.5, -9); 
    \node at (-3.2,-9) {19} ;
    \draw[thick] (-3, -10.5) -- (-1.5, -10.5); 
    \node at (-3.2,-10.5) {20} ;

    \draw[thick] (1.75, -4.5) -- (3.25, -4.5); 
    \node at (1.55,-4.5) {21} ;
    \draw[thick] (2, -6) -- (0.5, -6); 
    \node at (0.3,-6) {22} ;
    \draw[thick] (3, -6) -- (4.5, -6); 
    \node at (2.8,-6) {23} ;
    \draw[thick] (3, -7.5) -- (4.5, -7.5); 
    \node at (2.8,-7.5) {24} ;
    \draw[thick] (1.75, -9) -- (3.25, -9); 
    \node at (1.55,-9) {25} ;
    \draw[thick] (5.75, -9) -- (4.25, -9); 
    \node at (4.05,-9) {26} ;

    \node[draw, circle,scale=0.8] (u1) at (-0.5, 0.5) {U5};  
    \draw (u1) -- (-0.5, 0);

    \node[draw, circle,scale=0.8] (r4) at (-0.5, -1) {U6};  
    \draw (r4) -- (-0.5, -1.5);

    \node[draw, circle,scale=0.8] (r5) at (-0.5, -2.5) {R4};  
    \draw (r5) -- (-0.5, -3);

    \node[draw, circle,scale=0.8] (r6) at (-2.7, -4) {R5};  
    \draw (r6) -- (-2.7, -4.5);

    \node[draw, circle,scale=0.8] (r7) at (3, -4) {R7};  
    \draw (r7) -- (3, -4.5);

    \node[draw, circle,scale=0.8] (r8) at (-2.7, -5.5) {R6};  
    \draw (r8) -- (-2.7, -6);

    \node[draw, circle,scale=0.8] (r9) at (4.2, -5.5) {R8};  
    \draw (r9) -- (4.2, -6);

    \node[draw, circle,scale=0.8] (r10) at (4.2, -7) {R9};  
    \draw (r10) -- (4.2, -7.5);


    \draw[dashed, thick] (0, 0) -- (0,-1.5); 
    \draw[dashed, thick] (0, -3) -- (0,-1.5); 
    \draw[dashed, thick] (0, -3) -- (-2.25,-4.5); 
    \draw[dashed, thick] (0, -3) -- (2.25,-4.5);
    \draw[dashed, thick] (-2.25,-6) -- (-2.25,-4.5); 
    \draw[dashed, thick] (-2.25,-6) -- (-2.25,-7.5); 
    \draw[dashed, thick] (3.75,-6) -- (2.25,-4.5); 
    \draw[dashed, thick] (1,-6) -- (2.25,-4.5); 
    \draw[dashed, thick] (3.75,-6) -- (3.75,-7.5); 
    \draw[dashed, thick] (3.75,-7.5) -- (5.25,-9); 
    \draw[dashed, thick] (3.75,-7.5) -- (2.25,-9); 
    \draw[dashed, thick] (-2.25,-9) -- (-2.25,-7.5);
    \draw[dashed, thick] (-2.25,-9) -- (-2.25,-10.5);

    \draw[->,very thick] (0.5, 0) -- (0.5, 0.8); 
    \draw[->,very thick] (0.5, -1.5) -- (0.5, -0.7);
    \draw[->,very thick] (0.5, -3) -- (0.5, -2.2);

    \draw[->,very thick] (-1.8, -6) -- (-1.8, -5.2);
    \draw[->,very thick] (0.8, -6) -- (0.8, -5.2);
        
    \draw[->,very thick] (-2.7, -7.5) -- (-2.7, -6.7);
    \draw[->,very thick] (-2.7, -9) -- (-2.7, -8.2);
    \draw[->,very thick] (-2.7, -10.5) -- (-2.7, -9.7);

    \draw[->,very thick] (2, -9) -- (2, -8.2);
    \draw[->,very thick] (5.5, -9) -- (5.5, -8.2);

\end{tikzpicture}
}
    \caption{Shared topology of distribution networks. Each distribution system consists of two programmable generators (U5 and U6 for Distribution Netowrk 1), six non-programmable generators (R4, R5, R6, R7, R8, and R9 for Distribution Network 1), and ten loads. Of these, three are flexible loads, located at nodes 13, 14, and 15.}
    \label{fig:Dsystem}
\end{figure}

In the transmission subsystem, four programmable generators and five flexible loads are managed by three aggregators. Aggregator 1 controls programmable power plants U1 and U2, as well as the flexible load at node 4. Aggregator 2 manages programmable power plants U3 and U4 and the flexible load at node 6. Aggregator 3 controls the flexible loads at nodes 2, 3, and 5.
In each distribution system, two programmable power plants and three flexible loads are managed by two aggregators. 
The first aggregator controls one programmable generator and the major flexible load, while the second aggregator controls the second programmable generator and the remaining two flexible loads.

Table \ref{tab:gen_data} provides technical data for the programmable generators, indicating their reference system (transmission or distribution), their technology, their capacity, and their generation, upward and downward regulation costs.
\begin{table}[ht!]
\small
    \centering
    \begin{tabular} {lclrrrr}
    \hline
        $u \in \mathcal{U}$ & System & Technology & $G_u \ [MW]$ & $C_u \ [\frac{\text{\texteuro}}{MWh}]$  & $C_u^\text{UP} \ [\frac{\text{\texteuro}}{MWh}]$ & $C_u^\text{DW} \ [\frac{\text{\texteuro}}{MWh}]$ \\
        \hline
        U1 &  T & Coal & 500 & 88 & 132.0 & 44.0 \\
        U2 &  T & Hydro & 200 & 72 & 108.0 & 36.0 \\
        U3 &  T & CCS Gas & 300 & 91 & 136.5 & 45.5 \\
        U4 &  T & Gas & 500 & 71 & 106.5 & 35.5 \\ 
        U5 & D1 & CHP & 10 & 85 & 127.5 & 42.5 \\ 
        U6 & D1 & Gas Turbine & 5 & 80 & 120.0 & 40.0 \\
        U7 & D2 & Hydro & 5 & 75 & 112.5 & 37.5 \\
        U8 & D2 & CHP & 15 & 86 & 129.0 & 43.0 \\  
        U9 & D3 & Gas Turbine & 20 & 82 & 123.0 & 41 \\ 
        U10 & D3 & Hydro & 5 & 73 & 109.5 & 36.5 \\  
        \hline 
    \end{tabular}
    \caption{Data of programmable generators $u \in \mathcal{U}$. CHP stands for Combined Heat and Power plant}
    \label{tab:gen_data}
\end{table}
The ten programmable generating units are seven thermal power plants (U1, U3, U4, U5, U6, U8, and U9) and three hydro power plants (U2, U7, and U10). 
U1 is coal-fired, while the other thermal power plants are gas-fired. 
The programmable power plants have been assigned technology-specific production costs as reported in \cite{iea2020projected}. In contrast to other gas-fired plants, unit U3 is characterised by the installation of carbon capture and storage (CCS) technologies, making it the most expensive power plant in the system with a production cost of 91 \texteuro/MWh. On the other hand, U4 is the cheapest unit, with a production cost of 71 \texteuro/MWh, slightly lower than U2, U7 and U10, whose production costs reflect the opportunity costs of water resources. 
The regulations costs are defined from the production costs, assuming that upward (downward) adjustments have a 50\% higher (lower) cost. 

In terms of possible bidding strategies, in this study we consider for each programmable unit three price alternatives in the day-ahead market, with mark-ups of 10\%, 20\% and 30\%, and three price alternatives for upward (downward) regulations, with mark-ups (mark-downs) of 10\%, 30\%, and 50\% resulting in 27 possible strategies for each unit. The possible bid prices offered by the programmable units are summarized in Table \ref{tab:gen_prices}.

\begin{table}[ht!]
\small
    \centering
    \begin{tabular} {l| rrr| rrr| rrr}
    \hline
        \multirow{2}{*}{$u \in \mathcal{U}$} & \multicolumn{3}{c|}{Day-ahead bid} & \multicolumn{3}{c|}{Up-regulation bid} & \multicolumn{3}{c}{Down-regulation bid} \\
        & +10\% & +20\% & +30\% & +10\% & +30\% &  +50\% & $-10\%$ & $-20\%$ & $-50\% $\\
        \hline       
U1 & 96.80 & 105.60 & 114.40 & 145.20 & 171.60	& 198.00 & 39.60
& 30.80 & 22.00 \\
U2 & 79.20 & 86.40 & 93.60 & 118.80	& 140.40 & 162.00 &	32.40 & 25.20 &	18.00
 \\
U3 & 100.10 & 109.20 & 118.30 & 150.15 & 177.45 & 204.75 & 40.95 & 31.85 &	22.75
 \\
U4 & 78.10 & 85.20 & 92.30 & 117.15	& 138.45 & 159.75 &	31.95 &	24.85 &	17.75
 \\
U5 & 93.50	& 102.00 &	110.50 & 140.25	& 165.75 &	191.25 & 	38.25 &	29.75 &	21.25 \\
U6 & 88.00 & 96.00 & 104.00 & 132.00	& 156.00 &	180.00	& 36.00	& 28.00	& 20.00 \\
U7 & 82.50	& 90.00 &	97.50	& 123.75 &	146.25 & 168.75	& 33.75	& 26.25 &	18.75 \\
U8 & 94.60 & 103.20	& 111.80 &141.90 &167.70 &	193.50 &	38.70 &	30.10 &	21.50 \\
U9 & 90.20	& 82.00	& 106.60 &	135.30	& 159.90 &	184.50 &	36.90 &	28.70	& 20.50 \\
U10 & 80.30	& 73.00	& 94.90	& 120.45 &	142.35 &	164.25 &	32.85 & 25.55 &	18.25 \\
\hline
    \end{tabular}
    \caption{Bid price alternatives $[\frac{\text{\texteuro}}{MWh}]$ of programmable generators $u \in \mathcal{U}$.}
    \label{tab:gen_prices}
\end{table}


Regarding the demand side, 
we consider flexible 5 loads in transmission and the 3 major loads in each distribution subsystem. 
The flexible loads can only be curtailed in the services market and up to a maximum fraction $\delta_n$, which we assume to be 20\%.
Flexible loads present curtailment bids that compete in the services market with the up-regulation bids of programmable units. 
We consider three bidding strategies for each flexible load.
Table \ref{tab:load_data} provides data on the flexible loads ${\cal N}^F$, indicating their reference system (transmission or distribution), their load on the DAM and their possible bid prices on the ASM.
\begin{table} [ht!]
\small
    \centering
    \begin{tabular} {lcrrrr}
    \hline
         $n \in \mathcal{N}^F$&  System &  Load $[MW]$ & Price 1 $[\frac{\text{\texteuro}}{MWh}]$ & Price 2 $[\frac{\text{\texteuro}}{MWh}]$ & Price 3 [\texteuro/MWh]
\\
\hline
         N2&  T&  192.2 & 95.00 & 142.50	& 228.00 \\
         N3&  T&  219.2 & 98.79 & 148.18	& 237.09 \\
         N4&  T&  219.9 & 93.53 & 140.30	& 224.48 \\
         N5&  T&   69.5 & 95.74 &	143.61  & 229.77 \\
         N6&  T&  293.4 & 94.22 & 141.34	& 226.14 \\
         N13& D1& 12.86 & 97.47 &	146.21	& 233.94 \\
         N14& D1& 12.98 & 99.00 & 148.50	& 237.60 \\
         N15& D1&  3.35 & 96.35 &	144.53	& 231.25 \\
         N27& D2& 12.86 & 92.60 &	138.90	& 222.24 \\
         N28& D2& 12.98 & 94.05 &	141.08	& 225.72 \\
         N29& D2&  3.35 & 91.54 &	137.31	& 219.69 \\
         N41& D3& 12.86 & 97.57 &	146.36	& 234.18 \\
         N42& D3& 12.98 & 99.10 &	148.65	& 237.84 \\
         N43& D3&  3.35 &100.01 & 150.02	& 240.02 \\
 \hline
    \end{tabular}
    \caption{Data of system loads.}
    \label{tab:load_data}
\end{table}

Non-programmable generation in the case study can be categorized into two main types: photovoltaic (PV) units and wind turbines.
In the transmission system, two wind farms are connected to nodes 7 and 8, contributing 85 MWh and 95 MWh, respectively, to the DAM. 
Additionally, a PV plant connected to node 1 generates 65 MWh in the DAM.
Within the distribution networks, smaller-scale wind turbines and solar PV plants contribute to non-programmable generation. 
Specifically, the total non-programmable generation is 20 MWh in Distribution Network 1, 19 MWh in Distribution Network 2, and 31 MWh in Distribution Network 3.
In terms of imbalance, the deviation between the real-time and the day-ahead net load is assumed to be between $-129$ MWh and 129 MWh. This interval is sampled by 7 scenarios, as shown in Table \ref{tab:imbalance}.
%
%
\begin{table}[ht!]
\small
    \centering
    \begin{tabular}{l| rrrrrrr}
    \hline
        & $s_1$ & $s_2$ & $s_3$ & $s_4$ & $s_5$ & $s_6$ & $s_7$ \\
        \hline
        $\Delta_s$               &129 & 86 & 43 & 0 & $-43$ & $-86$ &$-129$ \\   
        \hline 
        $\Delta_s^{\cal T}$      & 99 & 66 & 33 & 0 & $-33$ & $-66$ & $-99$ \\
        $\Delta_s^{{\cal D}_1}$  &  9 &  6 &  3 & 0 &  $-3$ &  $-6$ &  $-9$ \\    
        $\Delta_s^{{\cal D}_2}$  &  6 &  4 &  2 & 0 &  $-2$ &  $-4$ &  $-6$ \\   
        $\Delta_s^{{\cal D}_3}$  & 15 & 10 &  5 & 0 &  $-5$ & $-10$ & $-15$ \\  \hline
    \end{tabular}
    \caption{Scenarios for the system imbalance $\Delta_s$ and its allocation between the transmission network and the three distribution systems.}
    \label{tab:imbalance}
\end{table}

\subsection{Results}
\label{sec:results}
All computational experiments were conducted on an ASUS laptop equipped with a 3 GHz Intel Core i7-5500U processor and 4 GB of RAM, using the GAMS 24.7.4 environment and the Gurobi solver.
Table \ref{tab:Model_Size} presents, for each coordination scheme, the dimensions of the optimization problem solved by each aggregator.
The table provides details on the number of constraints (\#Cons), binary variables (\#Bin), continuous variables (\#Var), and the average CPU time required at each iteration of the solution algorithm.
Note that aggregators 1 and 2 are grouped in the table, as they manage the same number of flexible resources (i.e., two programmable generators and one flexible load in the transmission system), resulting in optimization problems of identical dimensions. 
Similarly, aggregators 4, 6, and 8, which each control two flexible loads in the distribution subsystems, are grouped together. The same applies to aggregators 5, 7, and 9, which manage one programmable generator and one flexible load each.
As can be seen, the CPU times for each aggregator in schemes A and B differ only slightly. 
However, in Scheme C, the modeling of the sequential clearing between the distribution and transmission services markets leads to a significant increase in the size of the aggregators' models, which in turn results in a notable increase in solution times.
Regarding convergence, Scheme A requires 5 iterations to reach the equilibrium solution, while both schemes B and C reach convergence in 6 iterations. 
Consequently, the total solution times are 5275 seconds for Scheme A, 6804 seconds for Scheme B, and 10458 seconds for Scheme C.

\begin{table}[ht!]
\small
    \centering
    \begin{tabular}{l|l|rrrr}
    \hline
       & & \multicolumn{4}{c}{Aggregator} \\
         Scheme&  &  1, 2&  3&  4, 6, 8&  5, 7, 9\\
         \hline 
         A&  \#Cons & 127,460&  59,420&  93,440&  86,500\\
         A&  \#Bin &  10,047&  10,039&  10,043&  10,040\\
         A&  \#Var &  66,580&  50,580&  58,580&  56,270\\
         A&  CPU &  155&  91&  110&  108
\\
         \hline
         B&  \#Cons 
&  127,650&  59,630&  93,650&  86,710\\
         B&  \#Bin &  10,047&  10,039&  10,043&  10,040\\
         B&  \#Var &  66,790&  50,790&  58,790&  56,480\\
         B&  CPU &  158&  95&  120&  121
\\
         \hline
 C& \#Cons
& 490,651& 482,970& 488,458& 486,319\\
 C& 
\#Bin
& 10,047& 10,039& 10,043& 10,040\\
 C& \#Var
& 144,842& 144,321& 144,451& 144,405\\
         C&  CPU&  240&  156&  179&  190
\\
         \hline
    \end{tabular}
    \caption{Number of constraints, binary variables and continuous variables, and average CPU time (seconds) required to solve the decision problem of each aggregator under each coordination scheme.}
    \label{tab:Model_Size}
\end{table}

For each coordination scheme, Table \ref{tab:eq_down_prices} and Table \ref{tab:eq_prices} show the equilibrium bidding strategies of the flexible resources. 
In particular, Table \ref{tab:eq_down_prices} provides the optimal prices for the DAM and down-regulation bids of programmable power plants, while Table \ref{tab:eq_prices} shows the up-regulation bids submitted by programmable generators and the curtailment bids submitted by flexible loads in different coordination schemes. 
%
%
\begin{table}[ht!]
\small
    \centering
    \begin{tabular}
{l| rr| rrr| rrr } 
\hline 
\multirow{2}{*}{Resource} & \multicolumn{2}{c|}{Scheme A} & \multicolumn{3}{c|}{Scheme B} & \multicolumn{3}{c}{Scheme C} \\
   &    DAM  &    ASM &    DAM  &  ASM-D &  ASM-T &    DAM  &  ASM-D & ASM-T \\ 
   \hline
U1 &  96.80  & 39.60 &  96.80  &        & 39.60 &  96.80  &        & 39.60 \\
U2 &  93.60  & 32.40 &  93.60  &        & 32.40 &  93.60  &        & 32.40 \\
U3 & 100.10  & NP & 100.10  &            & NP & 100.10  &        & NP \\
U4 &  92.30  & 31.95 &  92.30  &        & 31.95 &  92.30  &        & 31.95  \\
U5 &  93.50  & 38.25 & 93.50  & 38.25  &        & 93.50  &  38.25  & 38.25\\
U6 &  88.00  & 36.00 &  88.00  & 36.00 &        &  88.00  & 36.00 & 36.00\\
U7 &  82.50  & 33.75 & 82.50  & 33.75   &       & 82.50  & 33.75 & 33.75  \\
U8 &  94.60  & 38.70 & 94.60  & 38.70  &        &  94.60  & 38.70 & 38.70 \\
U9 & 90.20  & 36.90 & 90.20  & 36.90  &        & 90.20 & 36.90 & 36.90\\
U10 &  80.30  & 32.85 & 80.30  & 32.55   &       & 80.30 & 32.85  & 32.85  \\
\hline
\end{tabular}
\caption{Equilibrium prices $[\frac{\text{\texteuro}}{MWh}]$ for DAM and down-regulation bids of programmable generators in different coordination schemes. NP stands for "Not Presented" and indicates a bid that is not submitted to the ASM as the generation unit is not dispatched in the DAM.}
\label{tab:eq_down_prices}
\end{table}

\begin{table}[ht!]
\small
\centering
\begin{tabular}
{l| r| rr| rr }
\hline
\multirow{2}{*}{Resource} &  \multicolumn{1}{c|}{Scheme A} & \multicolumn{2}{c|}{Scheme B} & \multicolumn{2}{c}{Scheme C} \\
   &    ASM &   ASM-D &  ASM-T &   ASM-D &  ASM-T \\ \hline
U1 & 145.20 &         & 145.20 &         & 145.20 \\
U3 & 150.15 &         & 150.15 &         & 150.15 \\
N2 & 142.50 &         & 142.50  &         & 142.50 \\ 
N3 & 148.10 &         & 148.10  &         & 148.10 \\
N4 & 140.30 &         & 140.30  &         & 140.30 \\
N5 & 229.76 &         & 229.76  &         & 229.76  \\
N6 & 141.34 &         & 141.34  &         & 141.34  \\
N13 &146.21 &   97.57 &         &  146.21 & 146.21  \\ 
N14 &148.50 &   99.10 &         &  148.50 & 148.50  \\
N15 &144.53 &  144.53 &         &  144.53 & 144.53  \\
N27 &138.90 &  138.90 &         &  138.90 & 138.90  \\
N28 &141.08 &  141.08 &         &  141.08 & 141.08  \\
N29 &137.31 &  137.31 &         &  137.31 & 137.31  \\
N41 &234.18 &   97.57 &         &  234.18 & 234.18  \\
N42 &237.84 &   99.10 &         &  237.84 & 237.84  \\
N43 &240.02 &  100.01 &         &  240.02 & 240.02 \\
\hline
\end{tabular}
\caption{Equilibrium prices $[\frac{\text{\texteuro}}{MWh}]$ for up-regulation bids of programmable generators and curtailment bids of flexible loads in different coordination schemes. 
The prices of up-regulation bids for units U2, U4, U5, U6, U7, U8, U9, and U10 are not provided as these generation units are fully dispatched in the DAM and thus cannot submit any up-regulation bid.}
\label{tab:eq_prices}
\end{table}

As observed in Table \ref{tab:eq_down_prices}, variations in the ancillary services market configuration do not affect the strategic behavior of operators in the DAM: the programmable power plants submit the same bid prices in all three schemes, resulting in the same dispatch. Specifically, to meet the day-ahead net load of 1019 MWh, transmission units U2 and U4 and all distribution generators are fully dispatched. The remaining 259 MWh are supplied by unit U1, setting the DAM price at 96.80 \texteuro/MWh.
However, this day-ahead production plan does not comply with network constraints, specifically violating the transit limits of the two transmission lines connecting nodes 1 and 6, and nodes 2 and 5, respectively. Transmission network congestion creates opportunities for market power exercise.
In particular, curtailing the load at node 5 becomes essential to alleviate the severe congestion of the transmission network, prompting the flexible load N5 to submit the maximum-priced curtailment bid (see Table \ref{tab:eq_prices}).

\begin{figure}
    \centering
    \begin{subfigure}{\textwidth}
        \centering
        \includegraphics[width=0.95\linewidth]{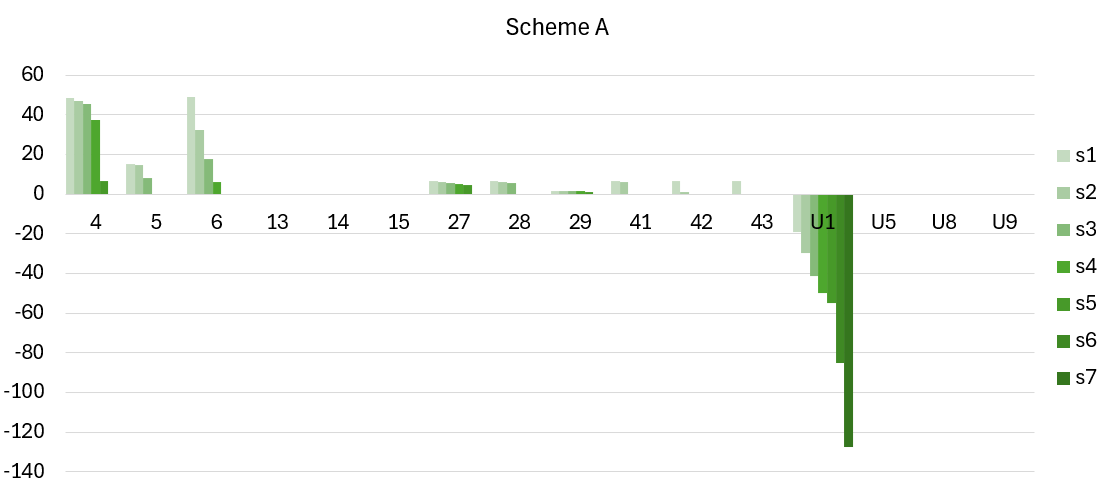}
        \caption{Dispatched quantities on the ASM in each imbalance scenario in coordination scheme A.}
        \label{fig:Q_A}
    \end{subfigure}

    \begin{subfigure}{\textwidth}
        \centering
        \includegraphics[width=0.95\linewidth]{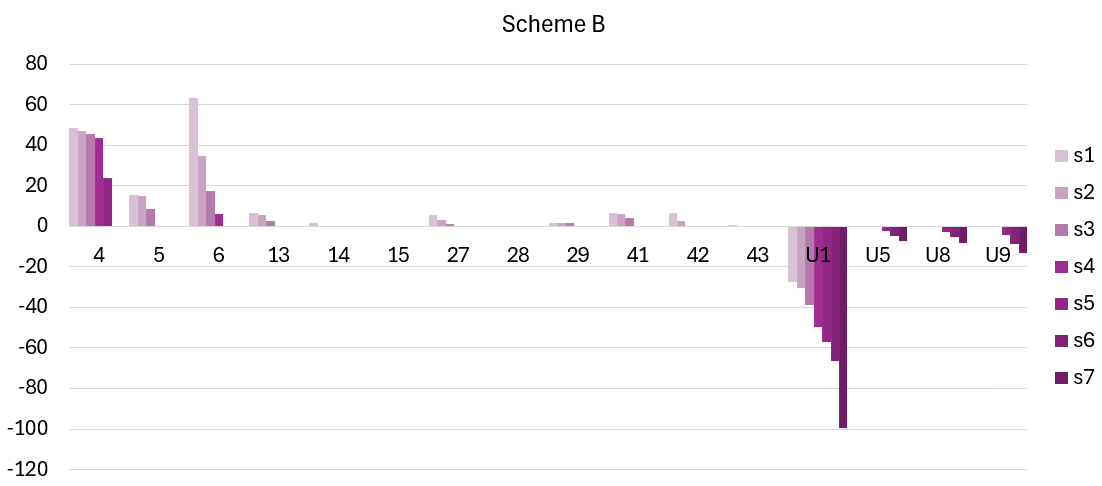}
        \caption{Dispatched quantities on the ASM in each imbalance scenario in coordination scheme B.}
        \label{fig:Q_B}
    \end{subfigure}

    \begin{subfigure}{\textwidth}
        \centering
        \includegraphics[width=0.95\linewidth]{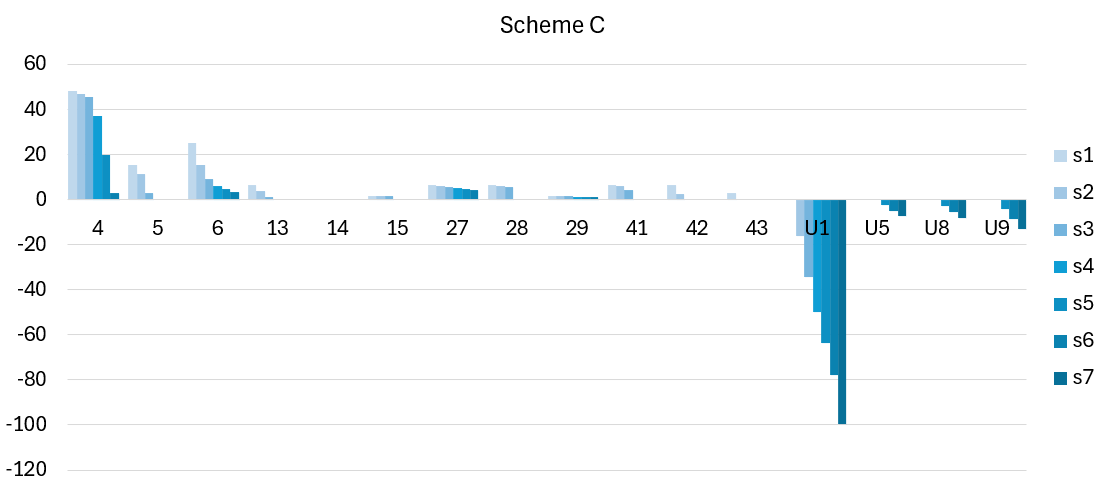}
        \caption{Dispatched quantities on the ASM in each imbalance scenario in coordination scheme A.}
        \label{fig:Q_C}
    \end{subfigure}
    
    \caption{Dispatched quantities on the ASM for flexible loads and programmable generators in each imbalance scenario under different coordination schemes.}
    \label{fig:Quantities}
\end{figure}

Regarding the distribution resources, different ASM configurations lead to varying strategic behaviors from market participants, with distribution flexible loads significantly increasing bid prices in schemes A and C.
Specifically, in coordination scheme A, the arbitrage opportunity created by the transmission network congestion extends to resources in distribution network 3, which submit maximum-priced bids in the common ancillary services market.
Conversely, in coordination scheme B, distribution resources can only participate in local markets aimed at solving network congestion and imbalances arising in the corresponding distribution system. As a result, the extremely high bid prices in the transmission system do not influence local competition among distribution flexible resources, which select the minimum-priced option that exceeds the clearing price $\lambda$.
Similar to scheme A, in scheme C, flexible loads in distribution system 3 also submit maximum-priced bids in both the distribution and transmission services markets, leveraging their critical role to alleviate congestion in the transmission network.
Fig. \ref{fig:Quantities} illustrates the curtailment and the regulations of programmable power plants on the ASM in each scenario under different coordination schemes.

Fig. \ref{fig:Q_A} illustrates that, under coordination scheme A, the ASM dispatch predominantly relies on resources connected to the transmission network.
In all scenarios, unit U1 is selected for downward regulation as it is the unit with the highest downward adjustment cost among those dispatched in the DAM.
In scenarios with negative imbalances (i.e., $s_5$, $s_6$, and $s_7$), reducing the power output of unit U1 effectively restores system balance. 
Conversely, in other scenarios, curtailing the production of U1 is necessary to alleviate congestion in the transmission network. 
This is further supported by load curtailments at nodes 4, 5, and 6, which help mitigating network congestion.
As seen in Fig. \ref{fig:Q_B}, the creation of dedicated services markets in scheme B favors the dispatch of distribution resources, which are employed to achieve local balancing within each distribution system.
Furthermore, the possibility for distribution resources to support the transmission system in scheme C increases the dispatch of distributed resources in the ASM, as shown in Fig. \ref{fig:Q_C}.

Finally, Table \ref{tab:costs} shows the total system costs for the three coordination schemes in each imbalance scenario.
It is observed that in scenarios $s_1$, $s_2$, and $s_3$, scheme B proves to be the most efficient architecture. In these cases, the significant system imbalance necessitates the acceptance of curtailment bids from resources located at node 5 of the transmission network, leading to a substantial increase in costs. Conversely, in scenarios with negative imbalances ($s_5$, $s_6$, and $s_7$), scheme A emerges as the more efficient option. The lower real-time load compared to the day-ahead forecast enables the alleviation of transmission network congestion without the need to accept curtailment offers at the aforementioned nodes. Overall, scheme B results in the lowest expected cost, totaling 6390.01 \texteuro. Scheme A is 5.13\% more expensive, while scheme C is 23.06\% more costly.

\begin{table}[ht!]
\small
    \centering
    \begin{tabular}{l| rrrrrrr |r}
    \hline
     Scheme   & $s_1$ & $s_2$ & $s_3$ & $s_4$ & $s_5$ & $s_6$ & $s_7$ & Expected Costs \\
     \hline
         A&  23137.87&  16934.08&  10822.99&  5003.63&  -445.83&  -3371.33&  -5057.00& 6717.77
\\
         B&  21227.42&  15752.79&  10357.93&  5014.85&  724.10&  -3338.82&  -5008.23& 6390.01
\\
         C&  26136.85&  18952.68&  11836.83&  5003.63&  1057.53&  -2888.56&  -5054.42& 7863.51
\\
\hline
    \end{tabular}
    \caption{Total system costs [\texteuro] for the provision of flexibility services in each scenario across the different coordination schemes.}
    \label{tab:costs}
\end{table}

\color{black}
\section{Conclusions}
\label{sec:concl}

Driven by the increasing demand for flexibility services due to the growing penetration of intermittent renewable generation, this paper explores the role of active distribution systems in providing ancillary services. Specifically, we investigate the potential of distributed energy resources to contribute to congestion management and balancing services under three distinct coordination schemes between services markets:
(A) a unified services market with coordinated management by TSO and DSOs;
(B) independent local markets for distribution resources, separate from the transmission services market;
(C) separate services markets for transmission and distribution, wherein unused distribution resources may participate in the transmission services market.

For each market configuration, we assess allocative efficiency and the potential for market power exploitation by modeling participants' bidding strategies in energy markets as a multi-leader common-follower game under uncertainty. In this framework, each market participant acts as a leader, optimizing the bidding strategies of its resources to maximize profits, while market operators act as followers, clearing energy markets to minimize total system costs. The resulting competitive game is expressed as a set of interdependent bilevel optimization problems, solved using an iterative procedure.

Two distinctive features of the proposed analysis make it particularly suited for assisting qualified operators in market design: (i) a two-stage stochastic programming modeling and (ii) a game-theoretical approach.
The two-stage stochastic programming framework provides a detailed representation of energy markets, capturing their sequential clearing stages and the dependencies between them. 
The game-theoretical approach allows evaluating the inherent vulnerability of electricity markets to market power exploitation, stemming from their oligopolistic structure and the sequential clearing mechanisms. 

The proposed methodology is validated through a case study, yielding several notable conclusions, some of which may initially seem counterintuitive. For example, theoretical expectations suggest that a single optimization process, as in Scheme A, would outperform separate processes like Schemes B and C. Moreover, Scheme C, where distribution resources support transmission services, is expected to be more cost-effective than Scheme B, where such support is absent. However, numerical results indicate that, in this case study, Scheme B incurs the lowest costs, followed by Scheme A and Scheme C.
The discrepancy arises from the market power potential induced by transmission network congestion. In such conditions, allowing distribution resources to participate in transmission service markets leads to higher bid prices from these resources. Conversely, when distribution resources are restricted from participating in transmission markets (as in Scheme B), the arbitrage opportunities stemming from transmission congestion do not influence bid prices in distribution markets.

The findings highlight the critical importance of accounting for market power dynamics in the design and coordination of energy markets. By offering a rigorous framework for evaluating market configurations, the proposed models provide valuable insights into improving economic efficiency and regulatory strategies in electricity markets.

\bigskip
\bigskip
\noindent
\textbf{Acknowledgments.}
This work has been financed by the Research Fund for the Italian Electrical System under the Three-Year Research Plan 2025-2027 (MASE, Decree n.388 of November 6th, 2024), in compliance with the Decree of April 12th, 2024.
\bigskip
\bibliographystyle{apacite}

\begin{thebibliography}{}

\bibitem [\protect \citeauthoryear {%
Allaz%
\ \BBA {} Vila%
}{%
Allaz%
\ \BBA {} Vila%
}{%
{\protect \APACyear {1993}}%
}]{%
allaz1993cournot}
\APACinsertmetastar {%
allaz1993cournot}%
\begin{APACrefauthors}%
Allaz, B.%
\BCBT {}\ \BBA {} Vila, J\BHBI L.%
\end{APACrefauthors}%
\unskip\
\newblock
\APACrefYearMonthDay{1993}{}{}.
\newblock
{\BBOQ}\APACrefatitle {Cournot competition, forward markets and efficiency} {Cournot competition, forward markets and efficiency}.{\BBCQ}
\newblock
\APACjournalVolNumPages{Journal of Economic theory}{59}{1}{1--16}.
\PrintBackRefs{\CurrentBib}

\bibitem [\protect \citeauthoryear {%
ARERA%
}{%
ARERA%
}{%
{\protect \APACyear {2017}}%
}]{%
ARERA2017}
\APACinsertmetastar {%
ARERA2017}%
\begin{APACrefauthors}%
ARERA.%
\end{APACrefauthors}%
\unskip\
\newblock
\APACrefYearMonthDay{2017}{}{}.
\newblock
\APACrefbtitle {{Delibera 05 maggio 2017 300/2017/R/eel. Prima apertura del mercato per il servizio di dispacciamento (MSD) alla domanda elettrica e alle unit\'a di produzione anche da fonti rinnovabili non gi\'a abilitate nonch\'e ai sistemi di accumulo. Istituzione di progetti pilota in vista della costituzione del testo integrato dispacciamento elettrico (TIDE) coerente con il balancing code europeo}.} {{Delibera 05 maggio 2017 300/2017/R/eel. Prima apertura del mercato per il servizio di dispacciamento (MSD) alla domanda elettrica e alle unit\'a di produzione anche da fonti rinnovabili non gi\'a abilitate nonch\'e ai sistemi di accumulo. Istituzione di progetti pilota in vista della costituzione del testo integrato dispacciamento elettrico (TIDE) coerente con il balancing code europeo}.}
\newblock
\begin{APACrefURL} \url{https://www.arera.it/atti-e-provvedimenti/dettaglio/17/300-17} \end{APACrefURL}
\PrintBackRefs{\CurrentBib}

\bibitem [\protect \citeauthoryear {%
Arosio%
\ \BBA {} Falabretti%
}{%
Arosio%
\ \BBA {} Falabretti%
}{%
{\protect \APACyear {2023}}%
}]{%
arosio2023participation}
\APACinsertmetastar {%
arosio2023participation}%
\begin{APACrefauthors}%
Arosio, M.%
\BCBT {}\ \BBA {} Falabretti, D.%
\end{APACrefauthors}%
\unskip\
\newblock
\APACrefYearMonthDay{2023}{}{}.
\newblock
{\BBOQ}\APACrefatitle {{DER Participation in Ancillary Services Market: An Analysis of Current Trends and Future Opportunities}} {{DER Participation in Ancillary Services Market: An Analysis of Current Trends and Future Opportunities}}.{\BBCQ}
\newblock
\APACjournalVolNumPages{Energies}{16}{5}{2443}.
\PrintBackRefs{\CurrentBib}

\bibitem [\protect \citeauthoryear {%
Dijk%
\ \BBA {} Willems%
}{%
Dijk%
\ \BBA {} Willems%
}{%
{\protect \APACyear {2011}}%
}]{%
dijk2011effect}
\APACinsertmetastar {%
dijk2011effect}%
\begin{APACrefauthors}%
Dijk, J.%
\BCBT {}\ \BBA {} Willems, B.%
\end{APACrefauthors}%
\unskip\
\newblock
\APACrefYearMonthDay{2011}{}{}.
\newblock
{\BBOQ}\APACrefatitle {The effect of counter-trading on competition in electricity markets} {The effect of counter-trading on competition in electricity markets}.{\BBCQ}
\newblock
\APACjournalVolNumPages{Energy Policy}{39}{3}{1764--1773}.
\PrintBackRefs{\CurrentBib}

\bibitem [\protect \citeauthoryear {%
Eid%
, Codani%
, Perez%
, Reneses%
\BCBL {}\ \BBA {} Hakvoort%
}{%
Eid%
\ \protect \BOthers {.}}{%
{\protect \APACyear {2016}}%
}]{%
eid2016managing}
\APACinsertmetastar {%
eid2016managing}%
\begin{APACrefauthors}%
Eid, C.%
, Codani, P.%
, Perez, Y.%
, Reneses, J.%
\BCBL {}\ \BBA {} Hakvoort, R.%
\end{APACrefauthors}%
\unskip\
\newblock
\APACrefYearMonthDay{2016}{}{}.
\newblock
{\BBOQ}\APACrefatitle {{Managing electric flexibility from Distributed Energy Resources: A review of incentives for market design}} {{Managing electric flexibility from Distributed Energy Resources: A review of incentives for market design}}.{\BBCQ}
\newblock
\APACjournalVolNumPages{Renewable and Sustainable Energy Reviews}{64}{}{237--247}.
\PrintBackRefs{\CurrentBib}

\bibitem [\protect \citeauthoryear {%
Garc{\i}%
, Conejo%
, Gabriel%
\BCBL {}\ \protect \BOthers {.}}{%
Garc{\i}%
\ \protect \BOthers {.}}{%
{\protect \APACyear {2006}}%
}]{%
garci2006electricity}
\APACinsertmetastar {%
garci2006electricity}%
\begin{APACrefauthors}%
Garc{\i}, R.%
, Conejo, A\BPBI J.%
, Gabriel, S.%
\BCBL {}\ \BOthersPeriod {.}\end{APACrefauthors}%
\unskip\
\newblock
\APACrefYearMonthDay{2006}{}{}.
\newblock
{\BBOQ}\APACrefatitle {Electricity market near-equilibrium under locational marginal pricing and minimum profit conditions} {Electricity market near-equilibrium under locational marginal pricing and minimum profit conditions}.{\BBCQ}
\newblock
\APACjournalVolNumPages{European Journal of Operational Research}{174}{1}{457--479}.
\PrintBackRefs{\CurrentBib}

\bibitem [\protect \citeauthoryear {%
Gulotta%
, Dacc{\`o}%
, Bosisio%
\BCBL {}\ \BBA {} Falabretti%
}{%
Gulotta%
\ \protect \BOthers {.}}{%
{\protect \APACyear {2023}}%
}]{%
gulotta2023opening}
\APACinsertmetastar {%
gulotta2023opening}%
\begin{APACrefauthors}%
Gulotta, F.%
, Dacc{\`o}, E.%
, Bosisio, A.%
\BCBL {}\ \BBA {} Falabretti, D.%
\end{APACrefauthors}%
\unskip\
\newblock
\APACrefYearMonthDay{2023}{}{}.
\newblock
{\BBOQ}\APACrefatitle {Opening of ancillary service markets to distributed energy resources: A review} {Opening of ancillary service markets to distributed energy resources: A review}.{\BBCQ}
\newblock
\APACjournalVolNumPages{Energies}{16}{6}{2814}.
\PrintBackRefs{\CurrentBib}

\bibitem [\protect \citeauthoryear {%
Hesamzadeh%
, Holmberg%
\BCBL {}\ \BBA {} Sarfati%
}{%
Hesamzadeh%
\ \protect \BOthers {.}}{%
{\protect \APACyear {2018}}%
}]{%
hesamzadeh2018simulation}
\APACinsertmetastar {%
hesamzadeh2018simulation}%
\begin{APACrefauthors}%
Hesamzadeh, M\BPBI R.%
, Holmberg, P.%
\BCBL {}\ \BBA {} Sarfati, M.%
\end{APACrefauthors}%
\unskip\
\newblock
\APACrefYearMonthDay{2018}{}{}.
\newblock
{\BBOQ}\APACrefatitle {Simulation and evaluation of zonal electricity market designs} {Simulation and evaluation of zonal electricity market designs}.{\BBCQ}
\newblock

\PrintBackRefs{\CurrentBib}

\bibitem [\protect \citeauthoryear {%
Hobbs%
\ \BBA {} Pang%
}{%
Hobbs%
\ \BBA {} Pang%
}{%
{\protect \APACyear {2007}}%
}]{%
hobbs2007nash}
\APACinsertmetastar {%
hobbs2007nash}%
\begin{APACrefauthors}%
Hobbs, B\BPBI F.%
\BCBT {}\ \BBA {} Pang, J\BHBI S.%
\end{APACrefauthors}%
\unskip\
\newblock
\APACrefYearMonthDay{2007}{}{}.
\newblock
{\BBOQ}\APACrefatitle {{Nash-Cournot equilibria in electric power markets with piecewise linear demand functions and joint constraints}} {{Nash-Cournot equilibria in electric power markets with piecewise linear demand functions and joint constraints}}.{\BBCQ}
\newblock
\APACjournalVolNumPages{Operations Research}{55}{1}{113--127}.
\PrintBackRefs{\CurrentBib}

\bibitem [\protect \citeauthoryear {%
Holmberg%
\ \BBA {} Lazarczyk%
}{%
Holmberg%
\ \BBA {} Lazarczyk%
}{%
{\protect \APACyear {2015}}%
}]{%
holmberg2015comparison}
\APACinsertmetastar {%
holmberg2015comparison}%
\begin{APACrefauthors}%
Holmberg, P.%
\BCBT {}\ \BBA {} Lazarczyk, E.%
\end{APACrefauthors}%
\unskip\
\newblock
\APACrefYearMonthDay{2015}{}{}.
\newblock
{\BBOQ}\APACrefatitle {{Comparison of congestion management techniques: Nodal, zonal and discriminatory pricing}} {{Comparison of congestion management techniques: Nodal, zonal and discriminatory pricing}}.{\BBCQ}
\newblock
\APACjournalVolNumPages{The Energy Journal}{}{}{145--166}.
\PrintBackRefs{\CurrentBib}

\bibitem [\protect \citeauthoryear {%
IEA%
}{%
IEA%
}{%
{\protect \APACyear {2020}}%
}]{%
iea2020projected}
\APACinsertmetastar {%
iea2020projected}%
\begin{APACrefauthors}%
IEA, N.%
\end{APACrefauthors}%
\unskip\
\newblock
\APACrefYearMonthDay{2020}{}{}.
\newblock
{\BBOQ}\APACrefatitle {Projected costs of generating electricity 2020} {Projected costs of generating electricity 2020}.{\BBCQ}.
\PrintBackRefs{\CurrentBib}

\bibitem [\protect \citeauthoryear {%
Losa%
\ \protect \BOthers {.}}{%
Losa%
\ \protect \BOthers {.}}{%
{\protect \APACyear {2021}}%
}]{%
losa2021platone}
\APACinsertmetastar {%
losa2021platone}%
\begin{APACrefauthors}%
Losa, I.%
, Monti, A.%
, Ginocchi, M.%
, Croce, V.%
, Bosco, F.%
, de Luca, E.%
\BDBL {}Petters, B.%
\end{APACrefauthors}%
\unskip\
\newblock
\APACrefYearMonthDay{2021}{}{}.
\newblock
{\BBOQ}\APACrefatitle {PLATONE: Towards a new open DSO platform for digital smart grid services and operation} {Platone: Towards a new open dso platform for digital smart grid services and operation}.{\BBCQ}
\newblock
\BIn{} \APACrefbtitle {CIRED 2021-The 26th International Conference and Exhibition on Electricity Distribution} {Cired 2021-the 26th international conference and exhibition on electricity distribution}\ (\BVOL\ 2021, \BPGS\ 2974--2978).
\PrintBackRefs{\CurrentBib}

\bibitem [\protect \citeauthoryear {%
McCormick%
}{%
McCormick%
}{%
{\protect \APACyear {1976}}%
}]{%
mccormick1976computability}
\APACinsertmetastar {%
mccormick1976computability}%
\begin{APACrefauthors}%
McCormick, G\BPBI P.%
\end{APACrefauthors}%
\unskip\
\newblock
\APACrefYearMonthDay{1976}{}{}.
\newblock
{\BBOQ}\APACrefatitle {{Computability of global solutions to factorable nonconvex programs: Part I -- Convex underestimating problems}} {{Computability of global solutions to factorable nonconvex programs: Part I -- Convex underestimating problems}}.{\BBCQ}
\newblock
\APACjournalVolNumPages{Mathematical programming}{10}{1}{147--175}.
\PrintBackRefs{\CurrentBib}

\bibitem [\protect \citeauthoryear {%
Migliavacca%
\ \protect \BOthers {.}}{%
Migliavacca%
\ \protect \BOthers {.}}{%
{\protect \APACyear {2021}}%
}]{%
migliavacca2021innovative}
\APACinsertmetastar {%
migliavacca2021innovative}%
\begin{APACrefauthors}%
Migliavacca, G.%
, Rossi, M.%
, Siface, D.%
, Marzoli, M.%
, Ergun, H.%
, Rodr{\'\i}guez-S{\'a}nchez, R.%
\BDBL {}others%
\end{APACrefauthors}%
\unskip\
\newblock
\APACrefYearMonthDay{2021}{}{}.
\newblock
{\BBOQ}\APACrefatitle {{The innovative FlexPlan grid-planning methodology: How storage and flexible resources could help in de-bottlenecking the European system}} {{The innovative FlexPlan grid-planning methodology: How storage and flexible resources could help in de-bottlenecking the European system}}.{\BBCQ}
\newblock
\APACjournalVolNumPages{Energies}{14}{4}{1194}.
\PrintBackRefs{\CurrentBib}

\bibitem [\protect \citeauthoryear {%
Migliavacca%
\ \protect \BOthers {.}}{%
Migliavacca%
\ \protect \BOthers {.}}{%
{\protect \APACyear {2017}}%
}]{%
migliavacca2017smartnet}
\APACinsertmetastar {%
migliavacca2017smartnet}%
\begin{APACrefauthors}%
Migliavacca, G.%
, Rossi, M.%
, Six, D.%
, D{\v{z}}amarija, M.%
, Horsmanheimo, S.%
, Madina, C.%
\BDBL {}others%
\end{APACrefauthors}%
\unskip\
\newblock
\APACrefYearMonthDay{2017}{}{}.
\newblock
{\BBOQ}\APACrefatitle {{SmartNet: H2020 project analysing TSO--DSO interaction to enable ancillary services provision from distribution networks}} {{SmartNet: H2020 project analysing TSO--DSO interaction to enable ancillary services provision from distribution networks}}.{\BBCQ}
\newblock
\APACjournalVolNumPages{{CIRED-Open Access Proceedings Journal}}{}{}{}.
\PrintBackRefs{\CurrentBib}

\bibitem [\protect \citeauthoryear {%
Parliament%
\ \BBA {} the Council of~the European~Union%
}{%
Parliament%
\ \BBA {} the Council of~the European~Union%
}{%
{\protect \APACyear {2019}}%
}]{%
EU}
\APACinsertmetastar {%
EU}%
\begin{APACrefauthors}%
Parliament, E.%
\BCBT {}\ \BBA {} the Council of~the European~Union.%
\end{APACrefauthors}%
\unskip\
\newblock
\APACrefYearMonthDay{2019}{}{}.
\newblock
\APACrefbtitle {{Directive (EU) 2019/944 of the European Parliament and of the Council of 5 June 2019 on Common Rules for the Internal Market for Electricity; Brussels, Belgium}.} {{Directive (EU) 2019/944 of the European Parliament and of the Council of 5 June 2019 on Common Rules for the Internal Market for Electricity; Brussels, Belgium}.}
\newblock
\begin{APACrefURL} \url{https://eur-lex.europa.eu/eli/dir/2019/944/oj} \end{APACrefURL}
\PrintBackRefs{\CurrentBib}

\bibitem [\protect \citeauthoryear {%
Rancilio%
, Rossi%
, Falabretti%
, Galliani%
\BCBL {}\ \BBA {} Merlo%
}{%
Rancilio%
\ \protect \BOthers {.}}{%
{\protect \APACyear {2022}}%
}]{%
rancilio2022ancillary}
\APACinsertmetastar {%
rancilio2022ancillary}%
\begin{APACrefauthors}%
Rancilio, G.%
, Rossi, A.%
, Falabretti, D.%
, Galliani, A.%
\BCBL {}\ \BBA {} Merlo, M.%
\end{APACrefauthors}%
\unskip\
\newblock
\APACrefYearMonthDay{2022}{}{}.
\newblock
{\BBOQ}\APACrefatitle {Ancillary services markets in europe: Evolution and regulatory trade-offs} {Ancillary services markets in europe: Evolution and regulatory trade-offs}.{\BBCQ}
\newblock
\APACjournalVolNumPages{Renewable and Sustainable Energy Reviews}{154}{}{111850}.
\PrintBackRefs{\CurrentBib}

\bibitem [\protect \citeauthoryear {%
Saha%
, Saleem%
\BCBL {}\ \BBA {} Roy%
}{%
Saha%
\ \protect \BOthers {.}}{%
{\protect \APACyear {2023}}%
}]{%
saha2023impact}
\APACinsertmetastar {%
saha2023impact}%
\begin{APACrefauthors}%
Saha, S.%
, Saleem, M.%
\BCBL {}\ \BBA {} Roy, T.%
\end{APACrefauthors}%
\unskip\
\newblock
\APACrefYearMonthDay{2023}{}{}.
\newblock
{\BBOQ}\APACrefatitle {Impact of high penetration of renewable energy sources on grid frequency behaviour} {Impact of high penetration of renewable energy sources on grid frequency behaviour}.{\BBCQ}
\newblock
\APACjournalVolNumPages{International Journal of Electrical Power \& Energy Systems}{145}{}{108701}.
\PrintBackRefs{\CurrentBib}

\bibitem [\protect \citeauthoryear {%
Sarfati%
, Hesamzadeh%
\BCBL {}\ \BBA {} Holmberg%
}{%
Sarfati%
\ \protect \BOthers {.}}{%
{\protect \APACyear {2018}}%
}]{%
sarfati2018increase}
\APACinsertmetastar {%
sarfati2018increase}%
\begin{APACrefauthors}%
Sarfati, M.%
, Hesamzadeh, M\BPBI R.%
\BCBL {}\ \BBA {} Holmberg, P.%
\end{APACrefauthors}%
\unskip\
\newblock
\APACrefYearMonthDay{2018}{}{}.
\newblock
{\BBOQ}\APACrefatitle {{Increase-Decrease Game under Imperfect Competition in Two-stage Zonal Power Markets--Part I: Concept Analysis}} {{Increase-Decrease Game under Imperfect Competition in Two-stage Zonal Power Markets--Part I: Concept Analysis}}.{\BBCQ}
\newblock

\PrintBackRefs{\CurrentBib}

\bibitem [\protect \citeauthoryear {%
Stanley%
, Johnston%
\BCBL {}\ \BBA {} Sioshansi%
}{%
Stanley%
\ \protect \BOthers {.}}{%
{\protect \APACyear {2019}}%
}]{%
stanley2019platforms}
\APACinsertmetastar {%
stanley2019platforms}%
\begin{APACrefauthors}%
Stanley, R.%
, Johnston, J.%
\BCBL {}\ \BBA {} Sioshansi, F.%
\end{APACrefauthors}%
\unskip\
\newblock
\APACrefYearMonthDay{2019}{}{}.
\newblock
{\BBOQ}\APACrefatitle {Platforms to support nonwire alternatives and DSO flexibility trading} {Platforms to support nonwire alternatives and dso flexibility trading}.{\BBCQ}
\newblock
\APACjournalVolNumPages{Consumer, Prosumer, Prosumager: How Service Innovations Will Disrupt the Utility Business Model}{}{}{111--126}.
\PrintBackRefs{\CurrentBib}

\bibitem [\protect \citeauthoryear {%
Strunz%
\ \protect \BOthers {.}}{%
Strunz%
\ \protect \BOthers {.}}{%
{\protect \APACyear {2014}}%
}]{%
strunztf}
\APACinsertmetastar {%
strunztf}%
\begin{APACrefauthors}%
Strunz, K.%
, Abbasi, E.%
, Fletcher, R.%
, Hatziargyriou, N.%
, Iravani, R.%
\BCBL {}\ \BBA {} Joos, G.%
\end{APACrefauthors}%
\unskip\
\newblock
\APACrefYearMonthDay{2014}{}{}.
\newblock
\APACrefbtitle {{TF C6. 04.02: TB 575—Benchmark Systems for Network Integration of Renewable and Distributed Energy Resources; CIGRE: Paris, France, 2014}} {{TF C6. 04.02: TB 575—Benchmark Systems for Network Integration of Renewable and Distributed Energy Resources; CIGRE: Paris, France, 2014}}\ \APACbVolEdTR{}{\BTR{}}.
\newblock
\APACaddressInstitution{}{ISBN 978-285-873-270-8}.
\PrintBackRefs{\CurrentBib}

\bibitem [\protect \citeauthoryear {%
Terna%
}{%
Terna%
}{%
{\protect \APACyear {2015}}%
}]{%
terna2012grid}
\APACinsertmetastar {%
terna2012grid}%
\begin{APACrefauthors}%
Terna, S.%
\end{APACrefauthors}%
\unskip\
\newblock
\APACrefYearMonthDay{2015}{}{}.
\newblock
\APACrefbtitle {Grid Code -- Chapter 4: Dispatching Regulations.} {Grid code -- chapter 4: Dispatching regulations.}
\newblock
\begin{APACrefURL} \url{https://www.terna.it/it/sistema-elettrico/ codici-rete/codice-rete-italiano} \end{APACrefURL}
\PrintBackRefs{\CurrentBib}

\bibitem [\protect \citeauthoryear {%
Utrilla%
, Davi-Arderius%
, Mart{\'\i}nez%
, Chaves-{\'A}vila%
\BCBL {}\ \BBA {} Arriola%
}{%
Utrilla%
\ \protect \BOthers {.}}{%
{\protect \APACyear {2020}}%
}]{%
utrilla2020large}
\APACinsertmetastar {%
utrilla2020large}%
\begin{APACrefauthors}%
Utrilla, F\BPBI D\BPBI M.%
, Davi-Arderius, D.%
, Mart{\'\i}nez, A\BPBI G.%
, Chaves-{\'A}vila, J\BPBI P.%
\BCBL {}\ \BBA {} Arriola, I\BPBI G.%
\end{APACrefauthors}%
\unskip\
\newblock
\APACrefYearMonthDay{2020}{}{}.
\newblock
{\BBOQ}\APACrefatitle {{Large-scale demonstration of TSO--DSO coordination: the CoordiNet Spanish approach}} {{Large-scale demonstration of TSO--DSO coordination: the CoordiNet Spanish approach}}.{\BBCQ}
\newblock
\BIn{} \APACrefbtitle {CIRED 2020 Berlin Workshop (CIRED 2020)} {Cired 2020 berlin workshop (cired 2020)}\ (\BVOL\ 2020, \BPGS\ 724--727).
\PrintBackRefs{\CurrentBib}

\bibitem [\protect \citeauthoryear {%
Zhang%
\ \BBA {} Kim%
}{%
Zhang%
\ \BBA {} Kim%
}{%
{\protect \APACyear {2010}}%
}]{%
zhang2010two}
\APACinsertmetastar {%
zhang2010two}%
\begin{APACrefauthors}%
Zhang, D.%
\BCBT {}\ \BBA {} Kim, S.%
\end{APACrefauthors}%
\unskip\
\newblock
\APACrefYearMonthDay{2010}{}{}.
\newblock
{\BBOQ}\APACrefatitle {A two stage stochastic equilibrium model for electricity markets with forward contracts} {A two stage stochastic equilibrium model for electricity markets with forward contracts}.{\BBCQ}
\newblock
\BIn{} \APACrefbtitle {2010 IEEE 11th International Conference on Probabilistic Methods Applied to Power Systems} {2010 ieee 11th international conference on probabilistic methods applied to power systems}\ (\BPGS\ 194--199).
\PrintBackRefs{\CurrentBib}

\bibitem [\protect \citeauthoryear {%
Zhang%
\ \BBA {} Xu%
}{%
Zhang%
\ \BBA {} Xu%
}{%
{\protect \APACyear {2013}}%
}]{%
zhang2013two}
\APACinsertmetastar {%
zhang2013two}%
\begin{APACrefauthors}%
Zhang, D.%
\BCBT {}\ \BBA {} Xu, H.%
\end{APACrefauthors}%
\unskip\
\newblock
\APACrefYearMonthDay{2013}{}{}.
\newblock
{\BBOQ}\APACrefatitle {Two-stage stochastic equilibrium problems with equilibrium constraints: modeling and numerical schemes} {Two-stage stochastic equilibrium problems with equilibrium constraints: modeling and numerical schemes}.{\BBCQ}
\newblock
\APACjournalVolNumPages{Optimization}{62}{12}{1627--1650}.
\PrintBackRefs{\CurrentBib}

\end{thebibliography}

\newpage
\appendix
\setcounter{table}{0}
\renewcommand{\thetable}{A.\arabic{table}}
\section{Reformulation of aggregators' problems as single-level optimization models}
\label{KKT}
This appendix provides the equivalent formulation of the bilevel problems introduced in Section \ref{sec:model} as single-level stochastic programs with complementarity constraints. 

\subsection{Coordination scheme A}
The decision problem of aggregator $i\in {\cal I}$ in coordination scheme A can be formulated as the following single-level optimization program with complementarity constraints:

\begin{subequations}
\label{MODEL:ASM_2_KKT}
{\allowdisplaybreaks
%
%
\begin{flalign}
\max & \quad \sum_{u \in {\cal U}_i} (\lambda - C_u) \ g_u + \sum_{s \in {\cal S}} \sigma_s \Bigg[
\sum_{u \in {\cal U}_i} \Bigg( \sum_{a \in {\cal A}^{\cal U, \uparrow}_u} B^{\cal U, \uparrow}_{u,a} \ x^{\cal U, \uparrow}_{u,a} 
- C^{\uparrow}_{u} \Bigg) \ g^{\uparrow}_{u,s} +  
\nonumber
\\ 
& + \sum_{n \in {\cal N}_i} \Bigg( \sum_{a \in {\cal A}^{\cal N, \downarrow}_n} B^{\cal N, \downarrow}_{n,a} \ x^{\cal N, \downarrow}_{n,a} - \lambda \Bigg) \ d^{\downarrow}_{n,s} + 
\sum_{u \in {\cal U}_i} \Bigg( C^{\downarrow}_u - \sum_{a \in {\cal A}^{\cal U, \downarrow}_u} B^{\cal U, \downarrow}_{u,a} \ x^{\cal U, \downarrow}_{u,a} \Bigg) \  g^{\downarrow}_{u,s}
\Bigg] \qquad \qquad
\end{flalign}
\begin{alignat}{2}
%
%
\text{s.t. } & x^{\cal U}_{u,a} \in \{0, 1\}, \quad a \in {\cal A}^{\cal U}_u, \ && u \in {\cal U}_i, 
\label{eq:A.1b}
\\
%
%
\quad & \sum_{a \in {\cal A}^{\cal U}_u} x^{\cal U}_{u,a} = 1, && u \in {\cal U}_i, 
\label{eq:A.1c}
\\
%
%
\quad & 0 \leq g_u \perp \Bigg( \nu_u - \lambda + \sum_{a \in {\cal A}^{\cal U}_u} B^{\cal U}_{u,a} \ x^{\cal U}_{u,a} \Bigg) \geq 0, && u \in {\cal U}_i,
\label{eq:A.1d}
\\
%
%
\quad & 0 \leq g_u \perp ( \nu_u - \lambda + b^{\cal U}_u ) \geq 0, && u \in {\cal U}_{-i}, 
\label{eq:A.1e}
\\
%
%
\quad & 0 \leq \nu_u \perp (G_u - g_u) \geq 0, && u \in {\cal U}, 
\label{eq:A.1f}
\\
%
%
\quad & \sum_{u \in {\cal U}} g_u =  \sum_{n \in {\cal N}} D_n - \sum_{r \in {\cal R}} W_r, 
\label{eq:A.1g}
\\
%
%
\quad & x^{\cal U, \uparrow}_{u,a} \in \{0, 1\}, && a \in {\cal A}^{\cal U, \uparrow}_u, u \in {\cal U}_i, 
\label{eq:A.1h}
\\
%
%
\quad & \sum_{a \in {\cal A}^{\cal U, \uparrow}_u} x^{\cal U, \uparrow}_{u,a} = 1, && u \in {\cal U}_i, 
\label{eq:A.1i}
\\
%
%
\quad & x^{\cal N, \downarrow}_{n,a} \in \{0, 1\}, && a \in {\cal A}^{\cal N, \downarrow}_n, n \in {\cal N}_i, 
\label{eq:A.1j}
\\
%
%
\quad & \sum_{a \in {\cal A}^{\cal N, \downarrow}_n} x^{\cal N, \downarrow}_{n,a} = 1, && n \in {\cal N}_i, 
\label{eq:A.1k}
\\
%
%
\quad & x^{\cal U, \downarrow}_{u,a} \in \{0, 1\}, && a \in {\cal A}^{\cal U, \downarrow}_u, u \in {\cal U}_i, 
\label{eq:A.1l}
\\
%
%
\quad & \sum_{a \in {\cal A}^{\cal U, \downarrow}_u} x^{\cal U, \downarrow}_{u,a} = 1, && u \in {\cal U}_i, 
\label{eq:A.1m}
\\
%
%
\quad & 0 \leq g^{\uparrow}_{u,s} \perp \ \Bigg( \beta_{u,s} - \alpha_s + \sum_{l \in {\cal L}} H_{l,n(u)} \ \mu_{l,s} + 
\sum_{a \in {\cal A}^{\cal U, \uparrow}_u} B^{\cal U, \uparrow}_{u,a} \ x^{\cal U, \uparrow}_{u,a}
\Bigg) \geq 0, && u \in {\cal U}_i, s \in {\cal S}, 
\label{eq:A.1n}
\\
%
%
\quad & 0 \leq g^{\uparrow}_{u,s} \perp \ \Bigg( \beta_{u,s} - \alpha_s + \sum_{l \in {\cal L}} H_{l,n(u)} \ \mu_{l,s} + 
b^{\cal U,\uparrow}_{u} \Bigg) \geq 0, && u \in {\cal U}_{-i}, s \in {\cal S}, 
\label{eq:A.1o}
\\
%
%
\quad & 0 \leq d^{\downarrow}_{n,s} \perp \Bigg( \gamma_{n,s} - \alpha_s + \sum_{l \in {\cal L}} H_{l,n} \ \mu_{l,s} + 
\sum_{a \in {\cal A}^{\cal N, \downarrow}_n} B^{\cal N, \downarrow}_{n,a} \ x^{\cal N, \downarrow}_{n,a}
\Bigg) \geq 0, && n \in {\cal N}_i, s \in {\cal S}, 
\label{eq:A.1p}
\\
%
%
\quad & 0 \leq d^{\downarrow}_{n,s} \perp \Bigg( \gamma_{n,s} - \alpha_s + \sum_{l \in {\cal L}} H_{l,n} \ \mu_{l,s} + 
b^{\cal N,\downarrow}_n \Bigg) \geq 0, && n \in {\cal N}_{-i}, s \in {\cal S}, 
\label{eq:A.1q}
\\
%
%
\quad & 0 \leq g^{\downarrow}_{u,s} \perp \Bigg( \phi_{u,s} + \alpha_s - \sum_{l \in {\cal L}} H_{l,n(u)} \ \mu_{l,s} - 
\sum_{a \in {\cal A}^{\cal U, \downarrow}_u} B^{\cal U, \downarrow}_{u,a} \ x^{\cal U, \downarrow}_{u,a}
\Bigg) \geq 0, \qquad && u \in {\cal U}_i, s \in {\cal S}, 
\label{eq:A.1r}
\\
%
%
\quad & 0 \leq g^{\downarrow}_{u,s} \perp \Bigg( \phi_{u,s} + \alpha_s - \sum_{l \in {\cal L}} H_{l,n(u)} \ \mu_{l,s} - 
b^{\cal U,\downarrow}_{u} \Bigg) \geq 0, && u \in {\cal U}_{-i}, s \in {\cal S}, 
\label{eq:A.1s}
\\
%
%
\quad & 0 \leq w^{\downarrow}_{r,s} \perp ( \chi_{r,s} +  \alpha
_s - \sum_{l \in {\cal L}} H_{l,n(r)} \ \mu_{l,s} ) \geq 0, && r \in {\cal R}, s \in {\cal S}, 
\label{eq:A.1t}
\\
%
%
\quad & 0 \leq \beta_{u,s} \perp (G_u - g_u - g^{\uparrow}_{u,s}) \geq 0, && u \in {\cal U}, s \in {\cal S}, 
\label{eq:A.1u}
\\
%
%
\quad & 0 \leq \gamma_{n,s} \perp (\delta_{n} \ \tilde{D}_{n,s} - d^{\downarrow}_{n,s}) \geq 0, && n \in {\cal N}, s \in {\cal S}, 
\label{eq:A.1v}
\\
%
%
\quad & 0 \leq \phi_{u,s} \perp (g_u - g^{\downarrow}_{u,s}) \geq 0, && u \in {\cal U}, s \in {\cal S}, 
\label{eq:A.1w}
\\
%
%
\quad & 0 \leq \chi_{r,s} \perp (\Tilde{W}_{r,s} - w^{\downarrow}_{r,s}) \geq 0, && r \in {\cal R}, s \in {\cal S}, 
\label{eq:A.1x}
\\
%
%
\quad & \sum_{u \in {\cal U}} g^{\uparrow}_{u,s} + \sum_{n \in {\cal N}} d^{\downarrow}_{n,s} 
- \sum_{u \in {\cal U}} g^{\downarrow}_{u,s}
- \sum_{r \in {\cal R}} w_{r,s}^{\downarrow}
= \Delta_{s}, && s \in {\cal S},
\label{eq:A.1y}
\\
%
%
\quad & 0 \leq \mu_{l,s} \perp 
\Bigg\{ \overline{F}_l - \sum_{n \in {\cal N}} H_{l,n} 
\Bigg[ \sum_{u \in {\cal U}_n} (g_u + g^{\uparrow}_{u,s} - g^{\downarrow}_{u,s})
+ &&
\nonumber \\
\quad & 
+ \sum_{r \in {\cal R}_n} (\Tilde{W}_{r,s} - w^{\downarrow}_{r,s})
 - ( \tilde{D}_{n,s} - d^{\downarrow}_{n,s}) \Bigg] 
 \Bigg\} \geq 0, \qquad
&&  l \in {\cal L}, s \in {\cal S}. 
\label{eq:A.1z}
\end{alignat}
%
%
}
\end{subequations}

In model \eqref{MODEL:ASM_2_KKT}, 
equations \eqref{eq:A.1d}$-$\eqref{eq:A.1g} are the KKT conditions of the DAM market.
Denoting with ${\cal U}_{-i} ={\cal U} \setminus{} {\cal U}_i$ the set of competitors' programmable power plants, in equation \eqref{eq:A.1e} the prices $b^{\cal U}_u$, $u \in {\cal U}_{-i}$, are the competitors' bids prices, which are fixed when solving the decision problem of Aggregator $i$.
Similarly, equations \eqref{eq:A.1n}$-$\eqref{eq:A.1z} are the KKT conditions of the ASM market.
We distinguish equations \eqref{eq:A.1n}, \eqref{eq:A.1p} and \eqref{eq:A.1r}, which include aggregator $i$'s prices, from equations \eqref{eq:A.1o}, \eqref{eq:A.1q} and \eqref{eq:A.1s}, including the competitors' bids prices $b^{{\cal U},\uparrow}_{u}, b^{{\cal U},\downarrow}_{u}$, $u \in {\cal U}_{-i}$, $b_n^{{\cal N},\downarrow}$, $n \in {\cal N}_{-i}$, which are fixed when solving the decision problem of Aggregator $i$.

\subsection{Coordination scheme B}
In this section, we provide the formulation of decision problem of Aggregator $i\in {\cal I}$ in coordination scheme B as a single-level optimization problem by deriving the optimality conditions of the distribution and transmission services market in coordination scheme B.
Specifically, the KKT conditions of the local services markets \eqref{eq:2e}$-$\eqref{eqn:2k} can be expressed as follows: for $s \in {\cal S}$ and $1 \leq k \leq K$
\begin{subequations}
\label{PROBLEM:B_D}
\allowdisplaybreaks {
\begin{alignat}{2}
& \sum_{u \in {\cal U}^{{\cal D}_k} } g^{\uparrow}_{u,s} 
+ \sum_{n \in {\cal N}^{{\cal D}_k} } d^{\downarrow}_{n,s} 
- \sum_{u \in {\cal U}^{{\cal D}_k} } g^{\downarrow}_{u,s} 
- \sum_{r \in {\cal R}^{{\cal D}_k} } w^{\downarrow}_{r,s}
= \Delta_{s}^{{\cal D}_k}  
= \qquad && \nonumber 
\\
& = \sum_{n \in {\cal N}^{{\cal D}_k} } ( \tilde{D}_{n,s} - D_n) 
- \sum_{r \in {\cal R}^{{\cal D}_k} } ( \tilde{W}_{r,s} - W_r), \qquad && 
\label{eqn:csB_17-D} 
\\
& 0 \leq \mu_{l,s} \perp 
\Bigg\{ \overline{F}_l - \sum_{n \in {\cal N}^{{\cal D}_k}} H_{l,n} \Bigg[ \sum_{u \in {\cal U}_n} (g_u + g^{\uparrow}_{u,s} - g^{\downarrow}_{u,s}) 
+ \qquad && \nonumber 
\\
& + \sum_{r \in {\cal R}_n} (\tilde{W}_{r,s} - w^{\downarrow}_{r,s}) - ( \tilde{D}_{n,s} - d^{\downarrow}_{n,s}) \Bigg] \Bigg\} \geq 0, \qquad && l \in {\cal L}^{{\cal D}_k},
\label{eqn:csB_18D} 
\\
%
& 0 \leq \beta_{u,s} \perp (G_u - g_u - g^{\uparrow}_{u,s}) \geq 0, \qquad && u \in {\cal U}^{{\cal D}_k},
\label{eqn:csB_20}
\\
& 0 \leq \phi_{u,s} \perp (g_u - g^{\downarrow}_{u,s}) \geq 0, \qquad && u \in {\cal U}^{{\cal D}_k},
\label{eqn:csB_22} \\
& 0 \leq \gamma_{n,s} \perp (\delta_{n} \ \tilde{D}_{n,s} - d^{\downarrow}_{n,s}) \geq 0, \qquad && n \in {\cal N}^{{\cal D}_k},
\label{eqn:csB_24}
\\
& 0 \leq \chi_{r,s} \perp (\tilde{W}_{r,s} - w^{\downarrow}_{r,s}) \geq 0, \qquad && r \in {\cal R}^{{\cal D}_k},
\label{eqn:csB_24RES}
\\
%
%
& 0 \leq g^{\uparrow}_{u,s} \perp \ \Bigg( \beta_{u,s} 
- \alpha^{\cal D}_{k,s}
+ \sum_{l \in {\cal L}^{{\cal D}_k}} H_{l,n(u)} \ \mu_{l,s} + 
\sum_{a \in {\cal A}^{\cal U, \uparrow}_u} B^{\cal U, \uparrow}_{u,a} \ x^{\cal U, \uparrow}_{u,a}
\Bigg) \geq 0, \qquad && u \in {\cal U}^{{\cal D}_k}_i,
\label{eqn:csB_19D}
\\
%
%
& 0 \leq g^{\uparrow}_{u,s} \perp \ \Bigg( \beta_{u,s} 
- \alpha^{\cal D}_{k,s}
+ \sum_{l \in {\cal L}^{{\cal D}_k}} H_{l,n(u)} \ \mu_{l,s} + 
b^{\cal U,\uparrow}_{u} \Bigg) \geq 0, \qquad && u \in {\cal U}^{{\cal D}_k}_{-i},
\label{eqn:csB_19DComp}
\\
%
& 0 \leq d^{\downarrow}_{n,s} \perp \Bigg( \gamma_{n,s} 
- \alpha^{\cal D}_{k,s} 
+ \sum_{l \in {\cal L}^{{\cal D}_k}} H_{l,n} \ \mu_{l,s} + 
\sum_{a \in {\cal A}^{\cal N, \downarrow}_n} B^{\cal N, \downarrow}_{n,a} \ x^{\cal N, \downarrow}_{n,a}
\Bigg) \geq 0, \qquad && n \in {\cal N}^{{\cal D}_k}_{i},
\label{eqn:csB_23D}
\\
%
& 0 \leq d^{\downarrow}_{n,s} \perp \Bigg( \gamma_{n,s} 
- \alpha^{\cal D}_{k,s} 
+ \sum_{l \in {\cal L}^{{\cal D}_k}} H_{l,n} \ \mu_{l,s} + 
b^{\cal N,\downarrow}_n \Bigg) \geq 0, \qquad && n \in {\cal N}^{{\cal D}_k}_{-i},
\label{eqn:csB_23DComp}
\\
& 0 \leq g^{\downarrow}_{u,s} \perp \Bigg( \phi_{u,s} + \alpha^{\cal D}
_{k,s} - \sum_{l \in {\cal L}^{\cal D}} H_{l,n(u)} \ \mu_{l,s} -  
\sum_{a \in {\cal A}^{\cal U, \downarrow}_u} B^{\cal U, \downarrow}_{u,a} \ x^{\cal U, \downarrow}_{u,a}
\Bigg) \geq 0, \qquad && u \in {\cal U}^{{\cal D}_k}_i, \label{eqn:csB_21D} \\
%
%
& 0 \leq g^{\downarrow}_{u,s} \perp \Bigg( \phi_{u,s} + \alpha^{\cal D}_{k,s} 
- \sum_{l \in {\cal L}^{{\cal D}_k}} H_{l,n(u)} \ \mu_{l,s} - 
b^{\cal U,\downarrow}_{u} \Bigg) \geq 0, \qquad && u \in {\cal U}^{{\cal D}_k}_{-i},
\label{eqn:csB_21DComp} \\
& 0 \leq w^{\downarrow}_{r,s} \perp \Bigg( \chi_{r,s} + \alpha^{\cal D}_{k,s} 
- \sum_{l \in {\cal L}^{{\cal D}_k}} H_{l,n(r)} \ \mu_{l,s} \Bigg) \geq 0, \qquad && r \in {\cal R}^{{\cal D}_k}.
\label{eqn:csB_21D}
\end{alignat}
}
\end{subequations}

Similarly, the KKT conditions of the transmission services market \eqref{eq:2l}$-$\eqref{eqn:2r} can be expressed as follows: for $s \in {\cal S}$
\begin{subequations}
\label{PROBLEM:B_T}
\allowdisplaybreaks {
\begin{alignat}{2}
& \sum_{u \in {\cal U^T}} g^{\uparrow}_{u,s} + \sum_{n \in {\cal N^T}} d^{\downarrow}_{n,s} 
- \sum_{u \in {\cal U^T}} g^{\downarrow}_{u,s}
- \sum_{r \in {\cal R}^{\cal T}} w^{\downarrow}_{r,s}
=  
\Delta_{s}^{\cal T}= 
\qquad && 
\nonumber 
\\ 
& = \sum_{n \in {\cal N}^{\cal T}} ( \tilde{D}_{n,s} - D_n) - \sum_{r \in {\cal R}^{\cal T}} ( \tilde{W}_{r,s} - W_r), \qquad &&   
\label{eqn:csB_17-T} 
\\
& 0 \leq \mu_{l,s} \perp 
\Bigg\{ \overline{F}_l -  \sum_{n \in {\cal N}} H_{l,n} \Bigg[ \sum_{u \in {\cal U}_n} (g_u + g^{\uparrow}_{u,s} - g^{\downarrow}_{u,s}) 
+ \qquad && 
\nonumber 
\\
& + \sum_{r \in {\cal R}_n} (\Tilde{W}_{r,s} - w^{\downarrow}_{r,s}) - ( \tilde{D}_{n,s} - d^{\downarrow}_{n,s}) \Bigg] \Bigg\} \geq 0, \qquad && l \in {\cal L}^{\cal T},  
\label{eqn:csB_18T} 
\\
& 0 \leq \beta_{u,s} \perp (G_u - g_u - g^{\uparrow}_{u,s})  \geq 0, \qquad && u \in {\cal U^T},
\label{eqn:csB_20}
\\
& 0 \leq \phi_{u,s} \perp (g_u - g^{\downarrow}_{u,s}) \geq 0, \qquad && u \in {\cal U^T},
\label{eqn:csB_22} 
\\
& 0 \leq \gamma_{n,s} \perp (\delta_{n} \ \tilde{D}_{n,s} - d^{\downarrow}_{n,s}) \geq 0, \qquad && n \in {\cal N^T}, 
\label{eqn:csB_24}
\\
& 0 \leq \chi_{r,s} \perp (\Tilde{W}_{r,s} - w^{\downarrow}_{r,s}) \geq 0, \qquad && r \in {\cal R^T}, \\
& 0 \leq g^{\uparrow}_{u,s} \perp \ \Bigg( \beta_{u,s} -  \alpha^{\cal T}
_s+ \sum_{l \in {\cal L}^{\cal T}} H_{l,n(u)} \ \mu_{l,s} + 
\sum_{a \in {\cal A}^{\cal U, \uparrow}_u} B^{\cal U, \uparrow}_{u,a} \ x^{\cal U, \uparrow}_{u,a}
\Bigg) \geq 0, \qquad && u \in {\cal U}^{\cal T}_i,
\label{eqn:csB_19T} \\
& 0 \leq g^{\uparrow}_{u,s} \perp \ \Bigg( \beta_{u,s} -  \alpha^{\cal T}
_s+ \sum_{l \in {\cal L}^{\cal T}} H_{l,n(u)} \ \mu_{l,s} +  b^{\cal U,\uparrow}_{u} \Bigg) \geq 0, \qquad && u \in {\cal U}^{\cal T}_{-i}, 
\label{eqn:csB_19TComp} 
\\
& 0 \leq d^{\downarrow}_{n,s} \perp \Bigg( \gamma_{n,s} - \alpha^{\cal T}
_s + \sum_{l \in {\cal L}^{\cal T}} H_{l,n} \ \mu_{l,s} +  
\sum_{a \in {\cal A}^{\cal N, \downarrow}_n} B^{\cal N, \downarrow}_{n,a} \ x^{\cal N, \downarrow}_{n,a}
\Bigg) \geq 0, \qquad && n \in {\cal N}^{\cal T}_i, 
\label{eqn:csB_23T}
\\
& 0 \leq d^{\downarrow}_{n,s} \perp \Bigg( \gamma_{n,s} - \alpha^{\cal T}
_s + \sum_{l \in {\cal L}^{\cal T}} H_{l,n} \ \mu_{l,s} +  b^{\cal N,\downarrow}_n \Bigg) \geq 0, \qquad && n \in {\cal N}^{\cal T}_{-i}, 
\label{eqn:csB_23T}
\\
& 0 \leq g^{\downarrow}_{u,s} \perp \Bigg( \phi_{u,s} + \alpha^{\cal T}
_s - \sum_{l \in {\cal L}^{\cal T}} H_{l,n(u)} \ \mu_{l,s} -  
\sum_{a \in {\cal A}^{\cal U, \downarrow}_u} B^{\cal U, \downarrow}_{u,a} \ x^{\cal U, \downarrow}_{u,a}
\Bigg) \geq 0, \qquad && u \in {\cal U}^{\cal T}_i, 
\label{eqn:csB_21T}
\\
& 0 \leq g^{\downarrow}_{u,s} \perp \Bigg( \phi_{u,s} + \alpha^{\cal T}
_s - \sum_{l \in {\cal L}^{\cal T}} H_{l,n(u)} \ \mu_{l,s} -  b^{\cal U,\downarrow}_{u} \Bigg) \geq 0, \qquad && u \in {\cal U}^{\cal T}_{-i}, \label{eqn:csB_21TComp} \\
& 0 \leq w^{\downarrow}_{r,s} \perp \Bigg( \chi_{r,s} + \alpha^{\cal T}
_s - \sum_{l \in {\cal L}^{\cal T}} H_{l,n(r)} \ \mu_{l,s} \Bigg) \geq 0, \qquad && r \in {\cal R}^{\cal T}. 
\label{eqn:csB_22T}
\end{alignat}
}
\end{subequations}

The bilevel program \eqref{MODEL:ASM_3.1} can be transformed into a single-level optimization problem by replacing the lower-level problems \eqref{eq:2e}$-$\eqref{eqn:2k} and  \eqref{eq:2l}$-$\eqref{eqn:2r} with equations 
\eqref{PROBLEM:B_D} and \eqref{PROBLEM:B_T}, respectively. 

\subsection{Coordination scheme C}
In this section, we provide the formulation of decision problem of aggregator $i\in {\cal I}$ in coordination scheme C as a single-level optimization problem by deriving the optimality conditions of the distribution and transmission services market in coordination scheme C.
Specifically, the KKT conditions of the local services markets \eqref{eq:3m}$-$\eqref{eqn:3s} can be expressed as follows: 
\begin{subequations}
\label{PROBLEM:KKT_C_1}
{\allowdisplaybreaks
\begin{alignat}{2}
& \sum_{u \in {\cal U}^{{\cal D}_k} } g^{\cal D,\uparrow}_{u,s} 
+ \sum_{n \in {\cal N}^{{\cal D}_k} } d^{\cal D,\downarrow}_{n,s} 
- \sum_{u \in {\cal U}^{{\cal D}_k} } g^{\cal D,\downarrow}_{u,s}
- \sum_{r \in {\cal R}^{{\cal D}_k} } w^{\cal D,\downarrow}_{r,s}
= \qquad && \nonumber \\
& = \sum_{n \in {\cal N}^{{\cal D}_k} } ( \tilde{D}_{n,s} - D_n) - \sum_{r \in {\cal R}^{{\cal D}_k} } ( \tilde{W}_{r,s} - W_r), \qquad &&  
\label{eqn:csC_17-D} 
\\
& 0 \leq \mu_{l,s}^{\cal D} \perp 
\Bigg\{ \overline{F}_l -  \sum_{n \in {\cal N}^{{\cal D}_k}} H_{l,n} \Bigg[ \sum_{u \in {\cal U}_n} (g_u + g^{\cal D,\uparrow}_{u,s} - g^{\cal D,\downarrow}_{u,s}) 
+ \qquad && 
\nonumber 
\\
& + \sum_{r \in {\cal R}_n} (\Tilde{W}_{r,s} - 
w^{\cal D,\downarrow}_{r,s}) - (\tilde{D}_{n,s} - d^{\cal D,\downarrow}_{n,s}) \Bigg] \Bigg\} \geq 0, \qquad && l \in {\cal L}^{{\cal D}_k}, \label{eqn:csC_18D} 
\\
& 0 \leq \beta^{\cal D}_{u,s} \perp (G_u - g_u - g^{\cal D,\uparrow}_{u,s}) \geq 0, \qquad && u \in {\cal U}^{{\cal D}_k},
\label{eqn:csC_20}
\\
& 0 \leq \gamma^{\cal D}_{n,s} \perp (\delta_{n} \ \tilde{D}_{n,s} - d^{\cal D,\downarrow}_{n,s}) \geq 0, \qquad && n \in {\cal N}^{{\cal D}_k},  
\label{eqn:csC_24}
\\
& 0 \leq \phi^{\cal D}_{u,s} \perp (g_u - g^{\cal D,\downarrow}_{u,s}) \geq 0, \qquad && u \in {\cal U}^{{\cal D}_k},
\label{eqn:csC_22} 
\\
& 0 \leq \chi^{\cal D}_{r,s} \perp (\Tilde{W}_{r,s} - w^{\cal D,\downarrow}_{r,s}) \geq 0, \qquad && r \in {\cal R}^{{\cal D}_k},  \\
& 0 \leq g^{\cal D,\uparrow}_{u,s} \perp \ \Bigg( \beta_{u,s}^{\cal D} 
- \alpha_{k,s}^{\cal D} 
+ \sum_{l \in {\cal L}^{{\cal D}_k}} H_{l,n(u)} \ \mu^{\cal D}_{l,s} + 
\sum_{a \in {\cal A}^{\cal U, \uparrow}_u} B^{\cal U, \uparrow}_{u,a} \ x^{\cal U, \cal D, \uparrow}_{u,a}
\Bigg) \geq 0, \qquad && u \in {\cal U}_i^{{\cal D}_k}, 
\label{eqn:csC_19D}
\\
& 0 \leq g^{\cal D,\uparrow}_{u,s} \perp \ \Bigg( \beta_{u,s}^{\cal D} 
- \alpha_{k,s}^{\cal D} 
+ \sum_{l \in {\cal L}^{{\cal D}_k}} H_{l,n(u)} \ \mu^{\cal D}_{l,s} +  b^{\cal U,\cal D,\uparrow}_{u} \Bigg) \geq 0, \qquad && u \in {\cal U}_{-i}^{{\cal D}_k}, 
\label{eqn:csC_19DComp}
\\
& 0 \leq d^{\cal D,\downarrow}_{n,s} \perp \Bigg( \gamma^{\cal D}_{n,s} 
- \alpha^{\cal D}_{k,s} 
+ \sum_{l \in {\cal L}^{{\cal D}_k}} H_{l,n} \ \mu^{\cal D}_{l,s} +  
\sum_{a \in {\cal A}^{\cal N, \downarrow}_n} B^{\cal N, \downarrow}_{n,a} \ x^{\cal N, \cal D, \downarrow}_{n,a}
\Bigg) \geq 0, \qquad && n \in {\cal N}^{{\cal D}_k}_{i},
\label{eqn:csC_23} \\
& 0 \leq d^{\cal D,\downarrow}_{n,s} \perp \Bigg( \gamma^{\cal D}_{n,s} 
- \alpha^{\cal D}_{k,s} 
+ \sum_{l \in {\cal L}^{{\cal D}_k}} H_{l,n} \ \mu^{\cal D}_{l,s} +  b^{\cal N,\cal D,\downarrow}_n \Bigg) \geq 0, \qquad && n \in {\cal N}^{{\cal D}_k}_{-i},
\label{eqn:csC_23Comp} 
\\
& 0 \leq g^{\cal D,\downarrow}_{u,s} \perp \Bigg( \phi^{\cal D}_{u,s} 
+ \alpha^{\cal D}_{k,s} 
- \sum_{l \in {\cal L}^{{\cal D}_k}} H_{l,n(u)} \ \mu^{\cal D}_{l,s} - 
\sum_{a \in {\cal A}^{\cal U, \downarrow}_u} B^{\cal U, \downarrow}_{u,a} \ x^{\cal U, \cal D, \downarrow}_{u,a}
\Bigg) \geq 0, \qquad && u \in {\cal U}^{{\cal D}_k}_i,
\label{eqn:csC_21} \\
& 0 \leq g^{\cal D,\downarrow}_{u,s} \perp \Bigg( \phi^{\cal D}_{u,s} 
+ \alpha^{\cal D}_{k,s} 
- \sum_{l \in {\cal L}^{{\cal D}_k}} H_{l,n(u)} \ \mu^{\cal D}_{l,s} - b^{\cal U,\cal D,\downarrow}_{u} \Bigg) \geq 0, \qquad && u \in {\cal U}^{{\cal D}_k}_{-i},
\label{eqn:csC_21Comp} \\
& 0 \leq w^{\cal D,\downarrow}_{r,s} \perp \Bigg( \chi^{\cal D}_{r,s} 
+ \alpha^{\cal D}_{k,s} 
- \sum_{l \in {\cal L}^{{\cal D}_k}} H_{l,n(r)} \ \mu^{\cal D}_{l,s} \Bigg) \geq 0, \qquad && r \in {\cal R}^{{\cal D}_k}.
\label{eqn:csC_21RES}
\end{alignat}
}\end{subequations}

Similarly, the KKT conditions of the transmission services market \eqref{eq:3ac}$-$\eqref{eqn:3am} can be expressed as follows: for $s \in {\cal S}$
\begin{subequations}
\label{PROBLEM:KKT_C_2}
{\allowdisplaybreaks
\begin{alignat}{2}
%
%
& \sum_{u \in {\cal U}} g^{\cal T,\uparrow}_{u,s} 
+ \sum_{n \in {\cal N}} d^{\cal T,\downarrow}_{n,s} 
- \sum_{u \in {\cal U}} g^{\cal T,\downarrow}_{u,s}
+ \sum_{r \in {\cal R}} w^{\cal T,\downarrow}_{r,s} 
= \Delta_{s}^{\cal T} = 
\qquad && \nonumber 
\\ 
&  = \sum_{n \in {\cal N}^{\cal T}} (\tilde{D}_{n,s} - D_n) 
   - \sum_{r \in {\cal R}^{\cal T}} (\tilde{W}_{r,s} - W_r), \qquad &&  
\label{eqn:csC_17-T} 
\\
%
%
& 0 \leq \mu_{l,s}^{\cal T} \perp \Bigg\{ \overline{F}_l - \sum_{n \in {\cal N^T}} H_{l,n} 
    \Bigg[ \sum_{u \in {\cal U}_n} (g_u + g^{\cal T,\uparrow}_{u,s} - g^{\cal T,\downarrow}_{u,s}) + 
    \sum_{r \in {\cal R}_n} (\tilde{W}_{r,s} + && \nonumber \\
& \qquad  - w^{\cal T}_{r,s}) - (\tilde{D}_{n,s} - d^{\cal T,\downarrow}_{n,s}) \Bigg] - \sum_{k=1}^K \sum_{n \in {\cal N}^{{\cal D}_k}} H_{l,n} \Bigg[ \sum_{u \in {\cal U}_n} (g_u + g^{\cal T,\uparrow}_{u,s}  + && \nonumber 
\\
& \qquad - g^{\cal T,\downarrow}_{u,s} + {g}^{\cal D,\uparrow}_{u,s} - {g}^{\cal D,\downarrow}_{u,s}) 
+ \sum_{r \in {\cal R}_n} (\tilde{W}_{r,s} - w^{\cal T}_{r,s} - w^{\cal D}_{r,s})
- ( \tilde{D}_{n,s} + && 
\nonumber 
\\
& \qquad - d^{\cal T,\downarrow}_{n,s} - {d}^{\cal D,\downarrow}_{n,s}) \Bigg] \Bigg\} \geq 0, \qquad  && l \in {\cal L^T},
\label{eqn:csC_18T} 
\\
%
%
& 0 \leq \beta^{\cal T}_{u,s} \perp (G_u - g_u - g^{\cal T,\uparrow}_{u,s}) \geq 0, \, \, && u \in {\cal U}^{\cal T},
\label{eqn:csC_20T}
\\
%
%
& 0 \leq \beta^{\cal T}_{u,s} \perp (G_u - g_u - g^{\cal T,\uparrow}_{u,s} - {g}^{\cal D,\uparrow}_{u,s}) \geq 0, \qquad && u \in {\cal U}^{{\cal D}_k},1 \leq k \leq K,
\label{eqn:csC_20TD} 
\\
%
%
& 0 \leq \gamma^{\cal T}_{n,s} \perp (\delta_{n} \ \tilde{D}_{n,s} - d^{\cal T,\downarrow}_{n,s}) \geq 0, \qquad && n \in {\cal N}^{\cal T},
\label{eqn:csC_24T}
\\
%
%
& 0 \leq \gamma^{\cal T}_{n,s} \perp (\delta_{n} \ \tilde{D}_{n,s} - d^{\cal T,\downarrow}_{n,s} - {d}^{\cal D,\downarrow}_{n,s}) \geq 0, \qquad && n \in {\cal N}^{{\cal D}_k},1 \leq k \leq K,
\label{eqn:csC_24TD}
\\
%
%
& 0 \leq \phi^{\cal T}_{u,s} \perp (g_u - g^{\cal T,\downarrow}_{u,s}) \geq 0, \qquad && u \in {\cal U}^{\cal T},
\label{eqn:csC_22T} 
\\
%
%
& 0 \leq \phi^{\cal T}_{u,s} \perp (g_u - g^{\cal T,\downarrow}_{u,s} - g^{\cal D,\downarrow}_{u,s}) \geq 0, \qquad && u \in {\cal U}^{{\cal D}_k}, 1 \leq k \leq K,
\label{eqn:csC_22TD} 
\\
%
%
& 0 \leq \chi^{\cal T}_{r,s} \perp (\tilde{W}_{r,s} - w^{\cal T}_{r,s}) \geq 0, \qquad && r \in {\cal R}^{{\cal T}},
\label{eqn:csC_22TRES} 
\\
%
%
& 0 \leq \chi^{\cal T}_{r,s} \perp (\tilde{W}_{r,s} - w^{\cal T}_{r,s} - w^{\cal D}_{r,s}) \geq 0, \qquad && r \in {\cal R}^{{\cal D}_k},1 \leq k \leq K,
\label{eqn:csC_22TDRES} 
\\
%
%
& 0 \leq g^{\cal T,\uparrow}_{u,s} \perp \ \Bigg( \beta_{u,s}^{\cal T} - \alpha_s^{\cal T} + \sum_{l \in {\cal L}^{\cal T}} H_{l,n(u)} \ \mu^{\cal T}_{l,s} +  \sum_{a \in {\cal A}^{\cal U, \uparrow}_u} B^{\cal U, \uparrow}_{u,a} \ x^{\cal U, \cal T,\uparrow}_{u,a}
\Bigg) \geq 0, \quad && u \in {\cal U}_i,
\label{eqn:csC_19T}
\\
%
%
& 0 \leq g^{\cal T,\uparrow}_{u,s} \perp \ \Bigg( \beta_{u,s}^{\cal T} - \alpha_s^{\cal T} + \sum_{l \in {\cal L}^{\cal T}} H_{l,n(u)} \ \mu^{\cal T}_{l,s} +  b^{\cal U,\cal T,\uparrow}_{u} \Bigg) \geq 0, \qquad && u \in {\cal U}_{-i},
\label{eqn:csC_19TComp}
\\
%
%
& 0 \leq d^{\cal T,\downarrow}_{n,s} \perp \Bigg( \gamma^{\cal T}_{n,s} - \alpha^{\cal T}_s + \sum_{l \in {\cal L}^{\cal T}} H_{l,n} \ \mu^{\cal T}_{l,s} + \sum_{a \in {\cal A}^{\cal N, \downarrow}_u} B^{\cal N, \downarrow}_{n,a} \ x^{\cal N, \cal T,\downarrow}_{n,a}
\Bigg) \geq 0,  && n \in {\cal N}_i,
\label{eqn:csC_23T}
\\
%
%
& 0 \leq d^{\cal T,\downarrow}_{n,s} \perp \Bigg( \gamma^{\cal T}_{n,s} - \alpha^{\cal T}_s + \sum_{l \in {\cal L}^{\cal T}} H_{l,n} \ \mu^{\cal T}_{l,s} + b^{\cal N,\cal T,\downarrow}_n \Bigg) \geq 0, \qquad && n \in {\cal N}_{-i},
\label{eqn:csC_23TComp}
\\
%
%
& 0 \leq {g}^{\cal T,\downarrow}_{u,s} \perp \Bigg( \phi^{\cal T}_{u,s} + \alpha^{\cal T}_s - \sum_{l \in {\cal L}^{\cal T}} H_{l,n(u)} \ \mu^{\cal T}_{l,s} - \sum_{a \in {\cal A}^{\cal U, \downarrow}_u} B^{\cal U, \downarrow}_{u,a} \ x^{\cal U, \cal T,\downarrow}_{u,a}
\Bigg) \geq 0,  && u \in {\cal U}_i, 
\label{eqn:csC_21TD} \\
%
%
& 0 \leq {g}^{\cal T,\downarrow}_{u,s} \perp \Bigg( \phi^{\cal T}_{u,s} + \alpha^{\cal T}_s - \sum_{l \in {\cal L}^{\cal T}} H_{l,n(u)} \ \mu^{\cal T}_{l,s} - {b}^{\cal U,\cal T,\downarrow}_{u} \Bigg) \geq 0, \qquad && u \in {\cal U}_{-i},
\label{eqn:csC_21TDComp} \\
%
%
& 0 \leq {w}^{\cal T,\downarrow}_{r,s} \perp \Bigg( \chi^{\cal T}_{r,s} + \alpha^{\cal T}_s - \sum_{l \in {\cal L}^{\cal T}} H_{l,n(r)} \ \mu^{\cal T}_{l,s}  \Bigg) \geq 0, \qquad && r \in {\cal R}.
\label{eqn:csC_21TDRES}
  \end{alignat}
}
\end{subequations}

The bilevel program \eqref{MODEL:ASM_3.2} can be transformed into a single-level optimization problem by replacing the lower-level problems \eqref{eq:3m}$-$\eqref{eqn:3s} and  \eqref{eq:3ac}$-$\eqref{eqn:3am} with equations 
\eqref{PROBLEM:KKT_C_1} and \eqref{PROBLEM:KKT_C_2}, respectively.

\end{document}